\documentclass[10pt,twoside]{article}
\usepackage{amsmath,amssymb}

\def\bf{\bfseries}

\newcommand{\traj}{{\rm Traj}}

\newcommand{\eps}{\varepsilon}

\newcommand{\dist}{{\rm dist}}

\newtheorem{theorem}{Theorem}
\newtheorem{itlemma}{Lemma}[section] 
\newtheorem{itproposition}[itlemma]{Proposition}
\newtheorem{itcorollary}[itlemma]{Corollary}
\newtheorem{itremark}[itlemma]{Remark}
\newtheorem{itdefinition}[itlemma]{Definition}
\newtheorem{itexample}[itlemma]{Example}

\newenvironment{lemma}{\begin{itlemma}\rm}{\end{itlemma}} 
\newenvironment{remark}{\begin{itremark}\rm}{\end{itremark}} 
\newenvironment{corollary}{\begin{itcorollary}\rm}{\end{itcorollary}}
\newenvironment{proposition}{\begin{itproposition}\rm}{\end{itproposition}}
\newenvironment{definition}{\begin{itdefinition}\rm}{\end{itdefinition}}
\newenvironment{example}{\begin{itexample}\rm}{\end{itexample}}

\newcommand{\bl}[1]{\begin{lemma}\label{#1}}
\newcommand{\br}[1]{\begin{remark}\label{#1}}
\newcommand{\bt}[1]{\begin{theorem}\label{#1}}
\newcommand{\bd}[1]{\begin{definition}\label{#1}}
\newcommand{\bp}[1]{\begin{proposition}\label{#1}}
\newcommand{\bc}[1]{\begin{corollary}\label{#1}}
\newcommand{\bfact}[1]{\begin{fact}\label{#1}}
\newcommand{\bex}[1]{\begin{example}\label{#1}}
\newcommand{\bem}[1]{\begin{example}\label{#1}}  
\newcommand{\ec}{\end{corollary}}
\newcommand{\eex}{\end{example}}
\newcommand{\eem}{\end{example}}
\newcommand{\el}{\end{lemma}}
\newcommand{\er}{\end{remark}}
\newcommand{\et}{\end{theorem}}
\newcommand{\ed}{\end{definition}}
\newcommand{\ep}{\end{proposition}}
\newcommand{\epr}{\end{proof}}
\newcommand{\bpr}{\begin{proof}}

\newcommand{\halmos}{\rule{1ex}{1.4ex}}

\newcommand{\beq}{\begin{eqnarray}}
\newcommand{\eeq}{\end{eqnarray}}
\newcommand{\beqn}{\begin{eqnarray*}}
\newcommand{\eeqn}{\end{eqnarray*}}
\newcommand{\bi}{\begin{itemize}}
\newcommand{\ei}{\end{itemize}}
\newcommand{\ben}{\begin{enumerate}}
\newcommand{\een}{\end{enumerate}}

\newenvironment{proof}{\noindent {\em Proof}.\ }{\hspace*{\fill}$\halmos$\medskip}
\newcommand{\twoif}[4]{
\left\{ \begin{array}{ll}#1&#2
\vspace{.4em}\\#3&#4\end{array}\right. }
\newcommand{\R}{{\mathbb R}}  
\newcommand{\N}{{\mathbb N}}  

\newcommand{\ds}{\displaystyle}

\newcommand{\cg}{{\cal G}}

\newcommand{\calo}{{\cal O}}
\newcommand{\kw}{\check{w}}
\newcommand{\dDelta}{\Delta\kern-.6em\Delta}
\newcommand{\ind}{1\kern-.4em 1}

\newcommand{\cvlg}{\check V_L[g]}

\makeatletter
\evensidemargin \oddsidemargin
\makeatother

\begin{document}
\thispagestyle{empty}
\setcounter{page}{1}

\noindent
{\footnotesize {\bf To appear in\\[-1mm]
{\em Dynamics of Continuous, Discrete and Impulsive Systems}}\\[-1.00mm]
http://monotone.uwaterloo.ca/$\sim$journal}
$~$ \\ [.3in]


\begin{center}
{\large\bf FURTHER RESULTS ON LYAPUNOV FUNCTIONS\\AND DOMAINS OF
ATTRACTION FOR PERTURBED ASYMPTOTICALLY STABLE SYSTEMS }

\vskip.20in

Michael Malisoff\\[2mm]
{\footnotesize Department of Mathematics, Louisiana State
University,
 Baton Rouge LA 70803-4918\\
$\mathtt{malisoff\, @\, lsu.edu}$,
$\mathtt{http://www.math.lsu.edu/\! \sim malisoff}$}
\end{center}

{\footnotesize \noindent {\bf Abstract.}  We present new theorems
characterizing robust Lyapunov functions and infinite horizon
value functions in optimal control as unique viscosity solutions
of partial differential equations.  We use these results to
further extend Zubov's method for representing domains of
attraction in terms of partial differential equation solutions.
 }
\vskip.2in

\section{Introduction}

\label{intro} The theories of Lyapunov functions and domains of
attraction form the basis for much of current work in stability
theory  (cf. \cite{BR01,  CGW01, GW00, S98}). An important result
in this area is the  {\em Zubov method} (cf. \cite{A83, H67,
KNB82, VV85, Z64}), which gives conditions under which the domain
of attraction of an asymptotically stable fixed point of $\dot
x=f(x)$ is $v^{-1}([0,1))$, where $v$ is the solution of the
 {\em Zubov equation}
\[
\ds Dv(x)\cdot f(x)=-H(x)[1-v(x)]\sqrt{1+||f(x)||^2},\; \; x\in
\R^N
\]
for suitable functions $H$. In \cite{CGW00a, CGW01, GW00}, Zubov's
method was extended to the important case of perturbed
asymptotically stable systems $\dot x=f(x,a)$ for which the fixed
point $0$ is stable under any perturbation $a$. These
perturbations, taken to be ${\cal A}:=\{ \text{measurable
functions}\,  \alpha:[0,\infty)\to A\}$ for a given compact set
$A$, are used to represent uncertainties  and exogenous effects,
rather than controls.  The main results in \cite{CGW00a, CGW01,
GW00} are partial differential equations (PDE) characterizations
for {\em robust domains of attraction} and {\em robust Lyapunov
functions} for perturbed dynamics (cf. $\S$\ref{defs} below for
the relevant definitions). Under the conditions in  \cite{CGW01},
the robust domain of attraction ${\cal D}_o$ for the perturbed
system $\dot x=f(x,a)$ is $v^{-1}([0,1))$, where $v$ is the unique
bounded viscosity solution of the {\em generalized Zubov equation}
\begin{equation}\label{gz}
  \ds \inf_{a \in A}\{-Dv(x)\cdot f(x,a)-g(x,a)+v(x)g(x,a)\}=0,\;
  \; x\in\R^N
\end{equation}
that vanishes at the origin, under certain restrictions on $g$.
The solution $v$ of (\ref{gz}) is a robust Lyapunov function for
the perturbed dynamics $f$. Also, if mild assumptions are added on
$g$ and $f$, then $v$ is locally Lipschitz (cf. $\S$\ref{lipsec}
below). The results from \cite{CGW01} form the basis for discrete
approximations of Lyapunov functions and ${\cal D}_o$ (cf.
\cite{CGW00a}).  In \cite{CGW00a, CGW01, GW00}, the function $g$
in (\ref{gz}) is assumed to satisfy\footnote{$B_r:=\{x\in\R^N:
||x||\le r\}$ for all $r>0$, where
$||x||=(x^2_1+x^2_2+\ldots+x^2_N)^{1/2}$.}
\begin{equation}\label{poss}
\begin{array}{l}
(i)\; \; \inf\{g(x,a): x\not\in B_r, a\in A\}\; \ge\;  g_o\; >\;
0\\ (ii)\; \;  g(x,a)>0\; \; \; \; \forall \; x\ne 0,\, \forall\;
a\in A\end{array}
\end{equation}
for some constants $g_o$ and $r$.
 However, one can easily find equations
(\ref{gz}) which admit several bounded solutions which are null at
the origin when (\ref{poss}) is not satisfied. Here is an example
where this occurs:

\bex{ex1}  Define the functions $f,g:\R\times[-1,1]\to\R$ by
\[
f(x,a):= \left\{
\begin{array}{ll}
1-\frac{a}{x}, &\text{if}\; \;  x\le -1,\\
-x+ax^2, &\text{if}\; \;  -1\le x\le 1,\\
\frac{a}{x}-1, &\text{if}\; \;  x\ge 1
\end{array}\right.\]
and
\[
g(x,a)\equiv
\gamma(x):=\twoif{|\sin(\pi x)|,}{-1\le x\le 1}{0,}{{\rm
otherwise}}.
\]
Then the robust domain of attraction ${\cal D}_o$ is $(-1,+1)$.
Let ${\cal M}$ denote the set of all measurable functions
$\alpha:[0,\infty)\to A:=[-1,+1]$ and $\traj_x(f)$ denote the set
of all solutions $\phi: [0,\infty)\to\R$ of $\dot y=f(y,a)$
starting at $x$ for $a\in {\cal M}$. One bounded
solution\footnote{We will understand our Zubov equation solutions
to be in the discontinuous viscosity sense (cf. $\S$\ref{defs},
Definition \ref{definitionvis}).} of the corresponding generalized
Zubov equation
\begin{equation}
\label{exeq1} \ds \inf_{a\in A}\{-f(x,a)\cdot Dv(x)\}
+(v(x)-1)\gamma(x)\; =\; 0
\end{equation}
on $\R$ is $v_1=1-\exp(-W_1)$, where $W_1$ is the  maximal cost
function
\[W_1(x)=\sup\left\{\int_0^\infty \gamma(\phi(t))\, {\rm d}t: \phi\in
\traj_x(f) \right\}\]  This follows from an elementary
consideration of semidifferentials.\footnote{By the arguments of
\cite{CGW01}, $v_1$ is a solution of (\ref{exeq1}) on any interval
of the form $(-\delta,+\delta)$ for $\delta\in(0,1)$, since on
this interval, $v_1$ is the Kru\v{z}kov transformation of a
maximal cost value function for a cost function satisfying
(\ref{poss}). Evidently $v_1$ is a solution of (\ref{exeq1})
outside $[-1,+1]$ where it is constant.  It remains to check the
semidifferential condition for $v_1$ to be  a solution of
(\ref{exeq1}) at $\pm 1$ (cf. \cite{BCD97}).  To check this
condition at $-1$, it suffices to check that (i) $-p-|p|\le 0$ for
all $p\in D^+(v_1)^\star(-1)$ and (ii) $-p-|p|\ge 0$ for all $p\in
D^-(v_1)_\star(-1)$. Condition (i) is trivial, while (ii) follows
because if $p\in D^-(v_1)_\star(-1)$, then
\[
\ds -p\ge \liminf_{h\downarrow
0}\frac{(v_1)_\star(-1+h)-(v_1)_\star(-1)-ph}{h}\ge 0,
\]
since $(v_1)_\star(-1+h)\le (v_1)_\star(-1)$ for $h>0$ near $0$. A
similar argument applies at $+1$.  This argument also shows that
the functions $v_2$ we are about to give are solutions of
(\ref{exeq1}) on $\R$.} Notice that $0\le v_1\le 1$ on $\R$ and
$v_1(0)=0$. Other bounded solutions of (\ref{exeq1}) on $\R$ are
\[
v_2(x)=\twoif{v_1(x),}{x>-1}{k,}{x\le -1},\; \; k>1\;
\text{constant}
\]
  As before,
$v_2(0)=0$, but $v_2(-1)>v_1(-1)$. Notice that $v_1$ and $v_2$ are
continuous at the origin. \eex

\br{rkintro}  Our hypotheses will imply  $v^{-1}([0,1))={\cal
D}_o$, where $v$ is the unique bounded  solution of the
corresponding generalized Zubov equation that satisfies $v(0)=0$
and that is continuous at the origin.   One hypothesis we will
make will be $g$ quasi-stability of $\R^N$, which is roughly the
condition that trajectories on $[0,\infty)$ with finite total cost
must approach the origin (cf. $\S$\ref{defs} for the  definition
of $g$ quasi-stability). In the previous example, this condition
is not satisfied, since $\phi(t)\equiv 1$ gives $0$ total costs.
On the other hand, all the other hypotheses of this note are
satisfied in this example, which shows that our quasi-stability
hypothesis cannot be dropped.\er

Example \ref{ex1} motivates our study of the solutions of the
generalized Zubov equation (\ref{gz}), and the sublevel sets of
these solutions, for general costs $g$. This note will develop
 {\em uniqueness} theory for solutions of (\ref{gz}), which includes
 the uniqueness characterizations given
in \cite{CGW01}, and also applies to cases which are not tractable
by the known results, e.g., cases where $g$ is non-Lipschitz or
violates (\ref{poss}).  We also study the {\em regularity} of
Zubov equation solutions,  and the sublevel sets of these
solutions,  under relaxed assumptions on $f$ and $g$. Our results
have the following  novel features:
 \bi
\item[1.]
Our results are based on   extensions of results in \cite{M00a,
M00, M01, MS02, MS00, S99} on uniqueness of solutions for the
infinite horizon Hamilton-Jacobi equation, namely, Theorem
\ref{thm1} and Propositions \ref{added1} and \ref{added2}  in
$\S$\ref{prelim}. The infinite horizon equation is the same as the
 exit time  equation. Whereas \cite{M01, MS02, MS00} assume
that the undiscounted infinite horizon Lagrangian is nonnegative,
our results apply for Lagrangians which could be null or negative.
It is natural to consider optimization with Lagrangians that take
both positive and negative values, to allow cost minimization  in
one part of the state space and  maximization  elsewhere. Results
on undiscounted exit problems with  negative Lagrangians were
given in \cite{M00}, which requires controllability to the
so-called positivity set of the Lagrangian. This controllability
condition is not needed below, nor do we need the uniform positive
lower bounds  on the interest rates used in \cite{S99a}. Theorem
\ref{thm1} does not put any  growth or lower bound assumptions on
solutions, nor does it require any controllability
 at the origin. On the other hand, the earlier
uniqueness characterizations for   Hamilton-Jacobi equations (cf.
\cite{BCD97, M00, M01, MS02, MS00, S99, S99a}) prove uniqueness of
solutions in classes of functions which are either proper or
bounded below by a finite constant, or which satisfy  an
asymptotic condition at the boundary of the domain. Therefore, in
addition to applying to problems with more general Lagrangians,
our results extend previous work by allowing more general
comparison functions, including functions which are negative and
neither bounded-from-above nor bounded-from-below.
\item[2.]
Our uniqueness theory for   (\ref{gz}) applies for general
functions $g$, and  gives stronger conclusions than the uniqueness
theory  of \cite{CGW00a, CGW01, GW00}.
 Clearly, a function $w$ is a solution of (\ref{gz}) exactly
when $-w$ is a solution of the usual infinite time Hamilton-Jacobi
equation (\ref{HJBE}) with interest rate $h=g$ and Lagrangian
$\ell=-g$. However, since the usual uniqueness results for
(\ref{HJBE}) require nonnegative $\ell$ or strictly positive $h$
(cf. \cite{BCD97, S99, S99a}), these results cannot in general be
applied to (\ref{gz}) when $g$ is nonnegative. Moreover, the known
results on (\ref{gz}) (cf. \cite{CGW00a, CGW01, GW00}) require
condition (\ref{poss}),  and a growth condition such as
Lipschitzness for $g$, and therefore cannot be applied to general
situations.  Since we allow general possibly non-Lipschitz costs
$g$, including cases where (\ref{poss}) is not satisfied (cf.
Remark \ref{compare} and $\S$\ref{illu}), our theory can be
regarded as an extension of the results of \cite{M00a} on optimal
control for non-Lipschitz dynamics. By allowing degenerate costs
$g$, we obtain solutions of (\ref{gz}) with  properties not found
in the Zubov equation solutions of \cite{CGW01} (cf.
$\S$\ref{lipsec} and Remark \ref{rknewly}). The unique solutions
of (\ref{gz}) are robust Lyapunov functions for $f$, and are the
Kru\v{z}kov transformations of  maximal cost type robust Lyapunov
functions  $V_L$ for $f$.  The functions $V_L$ are in turn  unique
solutions of
\begin{equation}\label{auxx}
\ds \inf_{a\in A}\{-f(x,a)\cdot Dv(x)-g(x,a)\}=0
\end{equation}
on ${\cal D}_o$. This generalizes the PDE characterizations for
(\ref{gz}) and (\ref{auxx})  in \cite{CGW01}. On the other hand,
the  uniqueness characterizations for (\ref{gz}) and (\ref{auxx})
in \cite{CGW01} all follow from the results given below.
Therefore, we obtain new classes of `flat' maximal cost type
robust Lyapunov functions, corresponding to degenerate cost
functions $g$,  which can be characterized as unique PDE solutions
(cf. Remark \ref{rknewly}). This leads to new characterizations of
${\cal D}_o$ as sublevel sets of value functions for degenerate
instantaneous costs (cf. Corollary \ref{prop2}).
 \ei
This note is organized as follows.  In $\S$\ref{defs}, we give
definitions and lemmas from \cite{BCD97, CGW01}. In
$\S$\ref{prelim}, we give  general uniqueness theory for solutions
of infinite horizon Hamilton-Jacobi equations. In $\S$\ref{appl},
we apply these results  to stabilization, and investigate the
regularity and uniqueness of  solutions of (\ref{gz}) and
(\ref{auxx}). These solutions are shown to be Lyapunov functions
with special properties not found in the Zubov equation solutions
in \cite{CGW01} (cf. $\S$\ref{lipsec}). Our results are based on
recent variations of the Filippov-Wa\v{z}ewski Relaxation Theorem,
which we review in Appendix A. In $\S$\ref{illu}, we illustrate
our results with an example from \cite{CGW01}. We close in
$\S$\ref{concl} by discussing extensions and directions for future
research.

\section{Definitions and Lemmas}
\label{defs} This note is concerned with perturbed systems of the
form $\dot x=f(x,a)$, $x(0)=x_o$ where the input $a$ represents
exogenous effects and uncertainties in the design of the control
system. As in  \cite{CGW01}, we will assume the following: \bi
\item[$(A_1)$]
$A$ is a nonempty compact metric space.
\item[$(A_2)$]
$f:\R^N\times A\to\R^N$ is bounded and continuous, where $N$ is a
fixed positive integer. Also, $f(x,a)$ is uniformly locally
Lipschitz, meaning, for each $R>0$, there is a constant $L_R>0$
such that if $||x||,||y||\le R$ and $a\in A$, then
$||f(x,a)-f(y,a)||\le L_R||x-y||$.
\item[$(A_3)$]
$f(0,a)=0$ for all $a\in A$. \ei In $\S$\ref{appl},  we  use
functions $g$ as the cost functions for {\em greatest cost}
Lyapunov functions, and we use cost functions $\ell$, which we
will refer to as {\bf Lagrangians},  when we study {\em least
cost} optimal control value functions  in $\S$\ref{prelim}. We
always assume that $g$ and $\ell$ satisfy the following: \bi
\item[$(A_4)$]
$g,\ell:\R^N\times A\to\R$ are continuous,  $\ell(0,a)=g(0,a)=0$
for all $a\in A$, and $g$ is nonnegative. \ei In our main
applications, we take $\ell(x,a)\le 0$ for all $(x,a)$, but our
theory applies to general continuous $\ell$ (cf. $\S$\ref{added}).
 Assumptions $(A_1)$-$(A_3)$ imply that for
each $a\in{\cal A}:=\{\, \text{measurable\  functions} \;  \alpha:
[0,\infty)\to A\}$  and $x_o\in\R^N$, the
 system \begin{equation}\label{eq}
 \dot x(t)=f(x(t),a(t)),\; \; \;   x(0)=x_o\end{equation} has a
unique (classical) solution defined on $[0,\infty)$, which will be
denoted by $\phi(\cdot,x_o,a)$, and called the {\bf trajectory} of
$f$  for $a$ starting at $x_o$. Elements of ${\cal A}$ are called
{\bf controls} or {\bf inputs}. For any $f$ and $A$ satisfying
$(A_1)$-$(A_3)$, ${\rm Traj}_{x_o}(f):=\{\phi(\cdot,x_o,a): a\in
{\cal A}\}$.

We  sometimes also use the  {\bf relaxed controls} ${\cal A}^r$,
which is the set of all measurable functions $\alpha:[0,\infty)\to
A^r$, where $A^r$ is the set of all Radon probability measures on
$A$ (cf. \cite{BCD97}). Recall (cf. \cite{BCD97}) that $A^r$ is a
compact metric space. By ${\cal A}^r\ni \alpha_n\to \alpha\in
{\cal A}^r$ weak-$\star$, we mean that for all $t\ge 0$ and  all
Lebesgue integrable functions $B:[0,t]\to C(A)$,
\begin{equation}
\label{neww} \ds \lim_{n\to\infty}\int_0^t\int_A (B(s))(a) \, {\rm
d}(\alpha_n(s))(a)\, {\rm d}s= \int_0^t\int_A (B(s))(a) \, {\rm
d}(\alpha(s))(a)\, {\rm d}s
\end{equation}
where $C(A)$ is the set of all real-valued continuous functions on
$A$.   For any $M\in \N$ and any function $\Phi:\R^N\times A\to
\R^M$, we define $\Phi^r:\R^N\times A^r\to\R^M$ by
$\Phi^r(x,m):=\int_A \Phi(x,a)\, {\rm d}m(a)$. Notice that
$(A_1)$-$(A_3)$ hold with $f$ replaced by $f^r$ and $A$ replaced
by $A^r$. We also use $\phi(\cdot,x,\alpha)$ to denote the
solution of $\dot y=f^r(y,\alpha)$ defined on $[0,+\infty)$
starting at $x$ for each $\alpha\in {\cal A}^r$. Therefore,
$\traj_x(f^r)=\{\phi(\cdot,x,\alpha):\alpha\in {\cal A}^r\}$. This
extends our original definition  of $\phi$, since we can view
${\cal A}$ as the subset of ${\cal A}^r$ consisting of all Dirac
probability measure valued relaxed controls. When $\alpha\in{\cal
A}^r$, we call $\phi(\cdot, x,\alpha)$ the {\bf relaxed
trajectory} of $f$ for $\alpha$ starting at $x$.  Recall the
following Compactness Lemma on relaxed controls (cf.
\cite{BCD97}): \bl{compactness} Let $(A_1)$-$(A_2)$ hold, let
$\{\alpha_n\}$ be a sequence in ${\cal A}^r$, and let $c>0$. Then
there exists a subsequence of $\{\alpha_n\}$ (which we do not
relabel) and an $\alpha\in {\cal A}^r$ such that (i)
$\alpha_n\to\alpha$ weak-$\star$ on $[0,c]$ and such that (ii) if
$x_n\to x$ in $\R^N$, then $\phi(\cdot, x_n, \alpha_n)\to
\phi(\cdot, x, \alpha)$ uniformly on $[0,c]$.\el

We sometimes also use a function $h$, representing a discount
rate, which we always assume satisfies \bi
\item[$(A_5)$]\ \
$h:\R^N\times A\to [0,+\infty)$ is continuous.  \ei

We say that $f$ is {\bf uniformly locally asymptotically stable
(ULAS)} provided that \bi\item[($\star\star$)] There are
$\beta_f\in {\cal KL}$ and  $r>0$ such that
$||\phi(t,x,\alpha)||\le \beta_f(||x||,t)$ for all  $x\in B_r$,
$t\ge 0$, and $\alpha\in {\cal A}$.\ei and that $f$ is {\bf
uniformly locally exponentially stable (ULES)} if ($\star\star$)
holds for $\beta_f(s,t)=Cse^{-\sigma t}$ for some positive
constants $C$ and $\sigma$.\footnote{Recall (cf. \cite{S98a}) that
${\cal K}^\infty$ is defined to be the set of
 all strictly increasing functions $F:[0,\infty)\to[0,\infty)$ which satisfy
 (i) $F(0)=0$ and (ii) $F(x)\to +\infty$
 as $x\to +\infty$.  Also, ${\cal KL}$ is the set of all
continuous functions
 $\beta:[0,\infty)\times [0,\infty)\to[0,\infty)$ for which
 (i) $\beta(\cdot, t)\in {\cal K}^\infty$ for all $t\ge 0$,
 (ii) $\beta(s,\cdot)$ is decreasing for all $s\ge 0$, and
(iii) $\beta(s,t)\to 0$ as $t\to +\infty$ for all $s\ge 0$.} When
$f$ is ULAS, we set $t(x,\alpha):=\inf\{t\ge 0:
||\phi(t,x,\alpha)||\le r\}$ for all $(x,\alpha)\in\R^N\times
{\cal A}^r$, where $\inf\emptyset$  is defined  to be $+\infty$.
Also,
\[
{\cal D}:=\left\{x\in\R^N: \lim_{t\to +\infty}\phi(t)= 0 \; \; \;
\forall \phi\in \traj_x(f)\right\}\] and \[{\cal
D}_o:=\left\{x\in\R^N: \sup_{\alpha\in{\cal A}
}t(x,\alpha)<+\infty\right\}.\] We call ${\cal D}$ the {\bf domain
of attraction} of $f$, and ${\cal D}_o$ is called the {\bf robust
domain of attraction} of $f$.   As shown in \cite{CGW01}, the sets
${\cal D}$ and ${\cal D}_o$ may differ for ULES dynamics
satisfying $(A_1)$-$(A_3)$. Since relaxed trajectories of $f$ can
be uniformly approximated by trajectories of $f$ on compact
intervals without changing the initial value, one shows that if
$f$ is ULAS, then so is $f^r$.  For each open set
$\cg\subseteq\R^N$, let
 $C^1(\cg)$ denote the set of all continuously
differentiable functions $F:\cg\to \R$. For each $S\subseteq\R^N$,
we set
\[w_\star(x)=\liminf_{S\ni y\to x}w(y)\; \; \; \; \text{and}\; \;
\; \; w^\star(x)=\limsup_{S\ni y\to x}w(y)\] for all $x\in S$ and
locally bounded functions $w:S\to\R$. Then $w_\star$ is lower
semicontinuous and  $w^\star$ is upper semicontinuous. Also,
$w_\star=w^\star=w$ at all points of continuity of $w$.
  We will find conditions under which the {\bf generalized Zubov
equation}
\begin{equation}\label{nHJBE}
  \ds \inf_{a \in A}\{-Dv(x)\cdot f(x,a)-g(x,a)+v(x)g(x,a)\}=0
\end{equation}
has a unique bounded  continuous solution $v$ on  $\R^N$ that
satisfies the two conditions $v(0)=0$ and $v^{-1}([0,1))={\cal
D}_o$. Also, we
 consider  solutions of
\begin{equation}
\label{HJBE} \ds \sup_{a\in A}\{-f(x,a)\cdot
Dv(x)-\ell(x,a)+h(x,a)v(x)\}=0
\end{equation}
which we refer to as the {\bf  infinite horizon Hamilton-Jacobi
equation} (for the dynamics $f$, Lagrangian $\ell$, and interest
rate $h$).   Since solutions of (\ref{nHJBE})-(\ref{HJBE}) may not
be differentiable or even continuous (cf. \cite{BCD97}), we
consider solutions of these equations in the viscosity sense, by
which we mean the following:

\bd{definitionvis}

Let $\cg\subseteq \R^N$ be open, $S\supseteq\cg$, $F:\R^N\times
\R\times \R^N\to\R$ be continuous,  and  $w: S\to\R$ be locally
bounded. We call $w$ a {\bf (discontinuous viscosity) solution} of
$F(x,w(x), Dw(x)) = 0$ on $\cg$ provided the following two
conditions hold:
\begin{itemize}\item[]\bi
\item[$(C_1)$]
If $\gamma\in C^1(\cg)$ and $x_o\in \cg$ is a local minimum of
$w_\star-\gamma$, then $F(x_o, w_\star(x_o), D \gamma(x_o))\,
\ge\, 0$.
\item[$(C_2)$]
If $\lambda\in C^1(\cg)$ and $x_1\in \cg$ is a local maximum of
$w^\star-\lambda$, then $F(x_1, w^\star(x_1), D \lambda(x_1))\,
\le\, 0$.
\end{itemize}\ei\ed
Let $\overline{\rm co}(S)$ denote the closed convex hull of any
set $S\subseteq\R^M$. We will set $\Phi(x,A) := \{\Phi(x,a): a\in
A\}$ for all $\Phi:\R^N\times A\to\R^M$, and ${\rm cl}(S)$ (resp.,
${\rm comp}(S)$) will denote the topological closure of $S$
(resp., $\R^N\setminus S$)  for each $S\subseteq\R^N$. For any
function $\Psi:\R^N\to[-\infty,+\infty]$, we define the {\bf
(effective) domain} of $\Psi$ to be ${\rm dom}(\Psi)=\{x\in\R^N:
\Psi(x)\in\R\}$, and we define ${\rm Trace}(\Psi):=\{\Psi(t): t\in
{\rm dom} (\Psi)\}$. Also, $\partial (B)$ denotes the boundary of
any set $B\subseteq\R^N$, and ${\rm Null}(F):=\{x: F(x)=0\}$ for
any real-valued  function $F$.   A set $\cg\subseteq\R^N$ will be
called {\bf (strongly) invariant} with respect to $f$ (or
$f$-invariant) provided the following holds: If $x\in\cg$, then
$\phi(t)\in\cg$ for all $t>0$ and $\phi\in {\rm Traj}_x(f)  $.  An
open  $f$-invariant set $\calo$ containing the origin is called
{\bf asymptotically null} ({\bf for $f$}) provided the following
holds: If $x\in \calo$ and $\phi\in\traj_x(f)$, then $\phi(t)\to
0$ as $t\to +\infty$ . One can show that ${\cal D}_o$ is
asymptotically null for $f$ when $f$ is ULES. Moreover, we have
the following, whose proof is in  \cite{CGW01}: \bl{lemma1} Let
$f$ be a ULAS dynamics satisfying $(A_1)$-$(A_3)$. Then: \bi
\item[a)]
${\rm cl} (B_r)\subseteq {\cal D}_o$, and ${\cal D}\subseteq {\rm
cl}({\cal D}_o)$.
\item[b)]
${\cal D}_o$ is open and asymptotically null for $f$.
\item[c)]
If $f(x,A)$ is convex for all $x\in \R^N$, then ${\cal D}={\cal
D}_o$. \ei\el

For each $x\in\R^N$, $\alpha\in {\cal A}^r$, and $h$ and $\ell$
satisfying $(A_4)$-$(A_5)$, we set
\[
\ds \tilde h(x,t,\alpha):=\int_0^t
h^r(\phi(s,x,\alpha),\alpha(s))\, {\rm d}s\] and
\[J[\ell,h](x,t,\alpha):=\int_0^t e^{-\tilde h(x,s,\alpha)}
\ell^r(\phi(s,x,\alpha),\alpha(s))\, {\rm d}s
\]
For cases where $h\equiv 0$, we use $J[\ell]$ to signify
$J[\ell,0]$. We also set \[\dist(p,S)=\inf\{||p-s||: s\in S\}\]
and $B_R(S) := \{x\in\R^N: \dist(x,S)\le R\}$ for each $R>0$,
$p\in \R^N$, and $S\subseteq \R^N$.  When $S=\{\bar x\}$, we write
$B_R(\bar x)$ instead of  $B_R(\{\bar x\})$.
 The following lemma follows from the proof
 of Theorem III.2.32 in \cite{BCD97}: \bl{lemma2} Let
$(A_1)$-$(A_5)$ hold, and $w$ be a solution of (\ref{HJBE}) on a
bounded open set $B\subseteq\R^N$.  Define  $\tau_q:{\cal A}\to
[0,\infty]$ and $T_\delta:\R^N\to [0,\infty]$ by
\begin{equation}
\label{ends}
\begin{array}{l}
\ds \tau_q(\beta)\; = \; \inf\{\, t\ge 0 :\; \phi(t,q, \beta)\,
\in\,
\partial B\}\\
T_\delta(p)\,  =\, \, \inf\{ t\ge 0: \phi(t,p,\alpha)\in
B_\delta(\partial B), \alpha\in {\cal A}\}\end{array}
\end{equation}
 for each $q\in B$ and $\delta>0$. Then: \bi
\item[(a)]
For all $q\in B$, $\beta\in {\cal A}$, and $r\in
[0,\tau_q(\beta))$,
 \[w^\star(q)\;  \le\;  \int_0^r e^{-\tilde
h(q,s,\beta)}\ell(\phi(s,q,\beta), \beta(s))\, {\rm d}s \;  +\;
e^{-\tilde h(q,r,\beta)}w^\star(\phi(r,q,\beta)).\]
\item[(b)]
For all  $q\in B$, $\delta\in(0,{\rm dist} (q,\partial B)/2)$, and
$t\in [0, T_\delta(q))$,
\begin{eqnarray}
w_\star(q)&\ge& \ds \inf_{\alpha\in {\cal A}}\left[\int_0^t
e^{-\tilde h(q,s,\alpha)}\ell(\phi(s,q,\alpha), \alpha(s))\, {\rm
d}s\right.\nonumber\\
 &+& \left. e^{-\tilde h(q,t,\alpha)} w_\star(\phi
(t,q,\alpha))\right].\nonumber\end{eqnarray} \ei\el A  function
$V:\calo\to\R$ on an open set $\calo\subseteq \R^N$ is called a
{\bf robust Lyapunov function} for $f$ provided (i) $V$ is
positive definite (meaning, $V(x)\ge 0$ for all $x$, and $V(x)=0$
iff $x=0$) and (ii) $V(x)> V(\phi(t,x,\alpha))$ for all $x\in
{\cal O}\setminus\{0\}$, $t>0$, and $\alpha\in {\cal A}$. An open
set $\calo\subseteq \R^N$ is called {\bf $g$ quasi-stable} (resp.,
{\bf relaxed $g$ quasi-stable}) provided the following condition
holds for each $x\in \calo$: If $\alpha\in {\cal A}$ and
$\int_0^\infty g(\phi(t,x,\alpha),\alpha(t))\, {\rm d}t<\infty$
(resp., $\alpha\in{\cal A}^r$ and $\int_0^\infty
g^r(\phi(t,x,\alpha),\alpha(t))\, {\rm d}t<\infty$), then
$\lim_{t\to+\infty}\phi(t,x,\alpha)=0$.

\br{complete} If we drop the assumption that $f$ is bounded but
keep the other hypotheses in $(A_1)$-$(A_3)$ the same, then
(\ref{eq}) has  a unique solution $\phi(\cdot,x_o,a)$ defined on a
maximal interval $[0,b)$, with $b>0$ depending on $x_o$ and $a\in
{\cal A}$. Recall that if $A$ is a compact metric space, then a
uniformly locally Lipschitz function $f:\R^N\times A\to\R^N$ (cf.
$(A_2)$) is called {\bf forward complete} provided that the
solution $\phi(\cdot, x_o,a)$ for (\ref{eq}) is defined on
$[0,\infty)$ for all $x_o\in \R^N$ and $a\in {\cal A}$ (cf.
\cite{LSW96}). Assumption $(A_2)$ implies that $f$ is forward
complete. As mentioned in \cite{GW00}, it is not very restrictive
to assume that $f$ is bounded, since we can replace $f$ by
$f/(1+||f||)$ without changing the trajectories. However,
normalizing in this way changes the design of the control system.
If we change $(A_1)$-$(A_3)$ by replacing the boundedness of $f$
with the assumption that $f$ is forward complete, then ${\cal
R}^T(S):=\{\phi(t,x,\alpha): 0\le t\le T, x\in S, \alpha\in {\cal
A}\}$ is bounded for each $T\in [0,\infty)$ and bounded set
$S\subseteq \R^N$ (cf. \cite{LSW96}). Using this fact, one can
show that all results in this note remain true if the boundedness
of $f$ is replaced by the forward completeness of $f$ (and the
other hypotheses are kept the same). The system $f$ will be
forward complete if, instead of assuming $f$ is bounded, we assume
$f(x,a)$ is globally Lipschitz in $x$ uniformly in $a$ (cf.
\cite{BCD97}).  \er

\br{con} It will not be necessary to assume that the sets $f(x,A)$
are convex.  Instead, we will instead apply Lemma \ref{lemma1}(c)
to the {\em relaxed} dynamics $f^r:\R^N\times A^r\to\R^N$, using
the fact that $f^r(x,A^r)\equiv \overline{{\rm co}}(f(x,A))$ is
convex for all $x\in\R^N$. In $\S$\ref{prelim}, we study solutions
of (\ref{HJBE}) on general asymptotically null sets, including
non-ULAS $f$  (cf. \cite{H67}, p. 191, for an example due to
Vinograd of asymptotically null sets for non-ULAS dynamics, with
no  controls, in which each trajectory asymptotically approaches
the origin).\er \br{fil} By the Filippov Selection Theorem,
relaxed $g$ quasi-stability and $g$ quasi-stability are equivalent
if $(f\times g) (x,A)$ is convex for each $x\in\R^N$ (cf.
\cite{BCD97}, $\S$VI.1). For sufficient conditions for $g$
quasi-stability for cascade systems, see \cite{MS02, MS00}. Notice
that both forms of $g$ quasi-stability hold for all choices of $g$
if ${\cal O}$ is asymptotically null for $f^r$. Also, both are
satisfied for $\calo=\R^N$ and ULES $f$ if (\ref{poss}) holds for
the constant $r$ in the ULES definition. Since ${\cal A}\subseteq
{\cal A}^r$, relaxed $g$ quasi-stability implies $g$
quasi-stability for any open set. In some applications, we also
assume $g$ is uniformly locally Lipschitz
 (cf. $(A_6)$ in $\S$\ref{appl} below), in which case $g$
quasi-stability and relaxed $g$ quasi-stability are equivalent
conditions for any open $\calo\subseteq\R^N$ (cf. Appendix A
below).\er

\section{Results on Infinite Horizon Problems}
\label{prelim} To develop  Lyapunov theory for ULAS dynamics, we
first study the infinite horizon value function
\begin{equation}
\label{inff} V_\infty(x):= \inf_{\alpha\in {\cal A}} J[\ell,
h](x,+\infty,\alpha)\; \; \in\; \; [-\infty,+\infty]
\end{equation} with the conventions that (i)
the infimum is only over those $\alpha$ for which
$J[\ell,h](x,+\infty,\alpha)$ converges in $\R\cup\{\pm \infty\}$
and (ii) $\inf\emptyset = +\infty$.  Our results on (\ref{inff})
will be of independent interest, because  we will allow general
continuous $\ell$, unlike the usual results, which  assume that
$\ell\ge 0$ everywhere, or that $h$ is uniformly bounded below by
a positive constant (cf. the remarks following the proof of
Theorem \ref{thm1} for a comparison of our results with the known
results). Since minimization of a function $F$ is equivalent to
maximizing $-F$, this allows problems where minimization takes
place in the part of the state space where $\ell\ge 0$, while
maximization takes place elsewhere. We will refer to the side
condition \bi\item[]\bi
\item[$(SC_w)$\ \ ]
$w(0)=0$, and $w$ is continuous at the origin \ei\ei In
particular, $(SC_w)$ implies that $w_\star(0)=w^\star(0)=0$, and
that for each $\eps>0$, there exists a $\delta>0$ for which
$w^\star(p)<\eps$ and $w_\star(p)>-\eps$ for all $p\in B_\delta$.
We will use this information in the proofs of Theorem \ref{thm1},
Theorem \ref{thm2}, and Proposition \ref{propii} below, to account
for the case where $w$ is discontinuous.

\subsection{Problems with Negative Lagrangians}
We begin with the case of {\em nonpositive} $\ell$, in which case
(\ref{inff}) involves the maximization of $J[-\ell,h]$. For the
extension to {\em general} Lagrangians, and to sets which are not
asymptotically null, see Remark \ref{negremark} and
$\S$\ref{added}.  Notice that $V_\infty$ might not be continuous
or even locally bounded, and so not a solution of (\ref{HJBE})
(cf. Remark \ref{negremark}, with $A=\{+1\}$). Nevertheless, we
have the following local comparison result, whose statement and
proof actually cover the cases where $\ell$ is everywhere
nonpositive or everywhere nonnegative:

\bt{thm1} Assume the following: \bi
\item[(1)]
$(A_1)$,$(A_2)$, $(A_4)$, $(A_5)$, ${\rm Null}(||f(0,\cdot)||)\ne
\emptyset$.
\item[(2)]
$\cg\subseteq\R^N$ is asymptotically null for $f$.
\item[(3)]
Either $\ell(x,a)\le 0$ for all $(x,a)\in \R^N\times A$, or else
$\ell(x,a)\ge 0$ for all $(x,a)\in \R^N\times A$.
\item[(4)]
 $w:\cg\to\R$ is
a solution of (\ref{HJBE}) on $\cg\setminus\{0\}$ that satisfies
$(SC_w)$.\ei Then $w\equiv V_\infty$ on $\cg$. \et \bpr We assume
$h\equiv 0$. The proof of the general case is similar to the one
we will now give. Let $\bar x\in \cg\setminus\{0\}$ and $\eps>0$
be given. It suffices to check that $w^\star(\bar x)\le
V_\infty(\bar x)$ and $w_\star(\bar x)\ge V_\infty(\bar x)$, since
$w_\star\le w^\star$ on $\cg$.  To do this, first note that
$J[\ell](\bar x,+\infty,\alpha)$ converges in $\R\cup\{\pm
\infty\}$ for all $\alpha\in {\cal A}$, by hypothesis (3). If
$\alpha\in {\cal A}$, then Lemma \ref{lemma2}, the asymptotic
nullness of $\cg$, and $(SC_w)$ give
\begin{eqnarray}
\ds w^\star(\bar x)&\le & \limsup_{t\to +\infty}\left[ \int_0^t
\ell(\phi(s,\bar x,\alpha),\alpha(s))\, {\rm
d}s+w^\star(\phi(t,\bar x,\alpha))\right] \nonumber\\&\le&
J[\ell](\bar x, +\infty,\alpha)+\limsup_{t\to +\infty}
w^\star(\phi(t,\bar x,\alpha))\; \; \le\; \; J[\ell](\bar x,
+\infty,\alpha) \nonumber
\end{eqnarray}
so $w^\star(\bar x)\le V_\infty(\bar x)$ follows by infimization.
Let $\eps>0$ be given.  It remains to construct $\hat\alpha\in
{\cal A}$ such that
\begin{equation}\label{goal}
  w_\star(\bar x)\; \; \ge\; \;  \int_0^M \ell(\phi(t,\bar x,\hat\alpha),
\hat\alpha(t))\,
  {\rm d}t\; +\; w_\star(\phi(M,\bar x,\hat\alpha))-\eps\; \; \; \; \forall
M\in\N
\end{equation}
Since $w_\star$ is lower semicontinuous and $\cg$ is
asymptotically null, a passage to the liminf as $M\to +\infty$ and
an infimization in (\ref{goal}) would then give $w_\star(\bar
x)\ge V_\infty(\bar x)-\eps$.  Since $\eps$ was arbitrary, we
would get  $w_\star(\bar x)\ge V_\infty(\bar x)$, as needed.  We
will now construct an $\hat\alpha$ satisfying (\ref{goal}) by
adapting ideas from \cite{MS02, MS00}.   We will assume that $w$
is continuous on $\cg$, the general case being exactly the same
but with $w_\star$ replacing $w$. Define $E_1, E_2, \ldots$ by
\begin{equation}
\label{ejay} E_j(t)\; \; \equiv\; \; \eps\, \left[\,
e^{-(j-1)}-e^{-(t+j-1)}\, \right] \; \; \text{for\ \  each}\; \; t
>0\end{equation}
By the choice of the $E_j$'s, \[E_1(1)+E_2(1)+\ldots +E_m(1)=
\eps(1-e^{-m})\] for all $m\in \N$. Set \[ {\cal Z}_1\; \; :=\; \;
\left\{\begin{array}{l}(t,\alpha)\in [0,1]\times {\cal A}:  \;
w(\bar x)\ge \int_0^t \ell(\phi(s,\bar x,\alpha), \alpha(s))\,
{\rm d}s\\
+w(\phi(t,\bar x,\alpha))- E_1(t)\end{array}\right\}
\]

The proof of the theorem in \cite{MS00}, which  we  review in
Appendix A, shows that ${\cal Z}_1$ contains an element of the
form $(1,\bar\alpha^{(1)})$. Set
\[c_1\equiv \bar \alpha^{(1)},\; \; \bar\gamma(0)=\bar x,\; \;
\bar\gamma(1)=\phi(1,\bar x,c_1).\] Next, we  will inductively
define the sets
\begin{equation}
\label{summa} {\cal Z}_{i+1}:=\left\{
\begin{array}{l}(t,\alpha)\in [0,1]\times
{\cal A}: \; w(\bar\gamma(i))-w(\phi(t,\bar\gamma(i),\alpha))\ge \\
\int_0^t \ell(\phi(s,\bar\gamma(i), \alpha),\alpha(s))\, {\rm d}s-
E_{i+1}(t)\end{array} \right\}\end{equation} Reapplying the
argument from \cite{MS00}, we inductively  obtain
\[(1,\bar\alpha^{(i)})\in {\cal Z}_i\; \; \text{and}\; \;
\bar\gamma(i+1):=\phi(1,\bar\gamma(i),\bar\alpha^{i+1})\in\cg\]
for all $i\in \N$. We inductively let $c_{i+1}$ be the
concatenation of $c_{i}\lceil[0,i]$ followed by
$\bar\alpha^{(i+1)}$.  Now sum both sides of the inequality in
(\ref{summa}), with the choices $\alpha=\bar\alpha^{i+1}$ and
$t=1$, over $i=0,1,2,\ldots, M-1$. Since \[ \sum_{i=0}^{
M-1}J[\ell](\bar\gamma(i),1,\bar\alpha^{i+1})= J[\ell](\bar
x,M,c_M)\] for all $M\in\N$, the choices of the $E_j$'s give
(\ref{goal}), with $\hat\alpha(t):=c_j(t)$ for $0\le t\le j$. \epr

\br{negremark} Notice that no growth conditions were required for
$\ell$, nor were there any requirements  on the {\em rate} the
trajectories approach the origin. In particular, we allow
non-Lipschitz $\ell$ and unstable $f$ (cf. Remark \ref{con}).
Condition (3) on $\ell$ in Theorem \ref{thm1} was used to force
the integral in (\ref{goal}) to converge in $\R\cup\{\pm
\infty\}$. However, for general $\ell$, $J[\ell](\bar
x,+\infty,\hat\alpha)$ may not converge in $\R$. This occurs if
for example
\[
\begin{array}{l}
N=1,\; \; \; 1\in A,\; \; \; f(x,1)\equiv -x^3,\\ h\equiv 0,\; \;
\; \ell(x,a)\equiv |x|, \; \; \; \hat\alpha\equiv 1\end{array}
\]
in which case
\[J[\ell](1,M,\hat\alpha)=\int_{0}^{M}(1+2t)^{-1/2}{\rm d}t,\]
 which does not converge in $\R$ as $M\to
+\infty$.  A different way to guarantee convergence of the costs
is to replace Condition (3) in Theorem \ref{thm1} with  one of the
following:\bi
\item[(A)]
$\ell$ is bounded, $\exists h_o>0$ s.t. $h(x,a)\ge h_o$ $\forall
(x,a)\in \R^N\times A$.
\item[(B)]
$|\ell(x,a)|\le h(x,a)$ $\forall (x,a)\in \R^N\times A$. \ei In
case (B), we  get \[J[|\ell|,h](x,+\infty,\alpha)\le
J[h,h](x,+\infty,\alpha)=1-\exp\left\{-J[h](x,+\infty,\alpha)\right\}\le
1\] for all $x\in\R^N$ and $\alpha\in {\cal A}$.   For case (A),
uniqueness characterizations for (\ref{HJBE}) are  known (cf.
\cite{BCD97}, Chapter III).  On the other hand, since we allow
general $h$ and $\ell$, case (B) is not covered by the usual
uniqueness  results. In both cases,
$J[\ell,h](x,+\infty,\alpha)\in \R$ for all $x\in \cg$ and
$\alpha\in {\cal A}$ and the argument we used to prove  Theorem
\ref{thm1} applies.  Further  extensions of Theorem \ref{thm1} are
given in $\S$\ref{added}. \er

\br{sw} If we add the assumption in Theorem \ref{thm1} that $f$ is
ULES, then the part of the proof of Theorem \ref{thm1} showing
that $w_\star(\bar x)\ge V_\infty(\bar x)$ can be simplified to
the following.  First assume $\ell(x,a)\le 0$ for all $x\in\R^N$
and $a\in A$. Let $r$ be as in the ULES definition, let $\eps>0$
be given, and choose $\delta\in (0,r)$ such that
$w_\star(p)>-\frac{\eps}{2}$ for all $p\in B_\delta$ (using  the
lower semicontinuity of $w_\star$ and $(SC_w)$).  For convenience,
we will assume that the constant $C$ in the ULES definition is
$1$, but a similar argument to the one we will now give applies if
$C>1$, by considering small overshoots.
 By
$(\star\star)$ (with $\beta_f(s,t)=se^{-\sigma t}$ for suitable
$\sigma>0$), $B_\delta$ is then $f$-invariant. Set
\[\tau(p,\mu,\alpha)=\inf\{t\ge 0: ||\phi(t,p,\alpha)||\le \mu\}\]
for all $p\in \cg$, $\mu>0$, and $\alpha\in {\cal A}$, and let
${\cal N}\subseteq\cg$ be an open neighborhood of $\bar x$. If
\[\bar T:= \sup\{\tau(\bar x,\delta,\alpha): \alpha\in {\cal
A}\}=+\infty,\] then by Lemma III.2 from \cite{SW96}, there exist
$q\in {\cal N}$ and $\beta\in {\cal A}$ such that
\[\tau(q,\delta/2,\beta)=+\infty.\] This contradicts the
asymptotically nullness of $\cg$.  It follows that $\bar
T<\infty$. By the convergence of all costs $J[\ell](x,+\infty,
\alpha)$, it suffices to show that
\[ {\cal Z}_{\bar T}\; \; :=\; \;
\left\{\begin{array}{l}(t,\alpha)\in [0,\bar T]\times {\cal A}:
w_\star(\bar x)\ge \int_0^t \ell(\phi(s,\bar x,\alpha),
\alpha(s))\, {\rm d}s\\+w_\star(\phi(t,\bar x,\alpha))-
\frac{\eps}{2}(1-e^{-t})\end{array}\right\}
\]
contains a pair $(\bar T,\bar\alpha)$, since then \[w_\star(\bar
x)\ge V_\infty(\bar x)-\eps.\]  The proof that ${\cal Z}_{\bar T}$
contains such a pair is similar to the argument of the theorem in
\cite{MS00} (cf. Appendix A below).  A similar simplification can
be made if, instead of assuming that $\ell$ is nonpositive, we
assume  $0\le \ell(x,a)\le ||x||$ for all $(x,a)$ (which implies
that $V_\infty$ is continuous at the origin). \er

\br{rk1} Whereas the results of \cite{MS00} assume the
controllability condition $STC\{0\}$ to the origin (cf.
\cite{BCD97}), the previous theorem does not require any
controllability.  For example, it applies to the ULES dynamics
\[f(x,a)= (-x_1,-x_2)+a(x^2_1, x^2_2),\;  A=[-1,+1],\] on $\cg
=(-1,+1)^2$,  where $STC\{0\}$ is not satisfied.  In the usual PDE
characterizations for (\ref{HJBE}) (cf. \cite{BCD97, M00a, M00,
MS00}), one also assumes that the comparison functions $w$ are
either bounded-from-below or that $w(x)\to +\infty$ as
$x\to\partial\cg$. However, these assumptions are not needed in
Theorem \ref{thm1}. Therefore, Theorem \ref{thm1} can be regarded
as an extension of the earlier results which takes a larger class
of possible dynamics and solutions into account.\er

\br{rk3}  Note that it was not necessary to assume any uniform
 bounds on the Lagrangians.  A uniform positive lower
bound  on  $\ell$ or $h$  is needed in the uniqueness
characterizations of \cite{BCD97, CGW01, S99}. Therefore, Theorem
\ref{thm1} also extends the earlier results  by allowing more
general Lagrangians and interest rates. Recall that $V_\infty$ may
not be locally bounded (cf. Remark \ref{negremark}, with
$A=\{+1\}$). However, the arguments of \cite{CGW01} establish that
if $(A_1)$-$(A_4)$ hold with $f$ ULES, $h\equiv 0$,  and $\ell=-g$
for $g$ uniformly locally Lipschitz (cf. $(A_6)$ in $\S$\ref{appl}
below), then $w=V_\infty$ is a solution of (\ref{HJBE}) on ${\cal
D}_o$ that satisfies $(SC_w)$. Since ${\cal D}_o$ is
asymptotically null for $f$, Theorem \ref{thm1} then gives a PDE
characterization for $V_\infty$ on ${\cal D}_o$.\er
\subsection{Problems with General Lagrangians}
\label{added} By combining the arguments of Theorem \ref{thm1}
with  \cite{MS00}, the following variant of Theorem \ref{thm1} is
easily shown:

\bp{added1}
 Assume the following:
\begin{itemize}
\item[(1)]
$(A_1)$, $(A_2)$, $(A_4)$, $h\equiv 0$, ${\rm
Null}(||f(0,\cdot)||)\ne \emptyset$.
\item[(2)]
$\ell$ is nonnegative, $\cg\subseteq\R^N$ is $\ell$ quasi-stable
and $f$-invariant.
\item[(3)]
$w:\cg\to\R$ is a bounded-from-below  solution of (\ref{HJBE}) on
$\cg\setminus\{0\}$ that satisfies $(SC_w)$.
\end{itemize}
Then $w\equiv V_\infty$ on $\cg$. \ep Proposition \ref{added1}
applies to the Fuller Example data from \cite{M00a, M00} on
$\cg=\R^2$ (in which $N=2$, $A=[-1,+1]$, $h\equiv 0$,
$f(x,a)=(x_2,a)$, $\ell(x,a)\equiv |x_1|^\gamma$, $\gamma>1$),
which is not tractable by means of the usual results on
first-order viscosity solutions. For a discussion of this
application, see \cite{M00a, MS02, MS00}. (For a different
uniqueness characterization that applies to the Fuller Example, as
well as  to other examples in which the quasi-stability condition
in (2) of Proposition \ref{added1} is not satisfied, see
\cite{M01}.)

The nonpositivity of $\ell$ in Theorem \ref{thm1} was used to
guarantee that the total cost $J[\ell,h](\bar
x,+\infty,\hat\alpha)$ for the constructed input $\hat\alpha$
converged in $\R\cup\{\pm \infty\}$.  This convergence is implied
by (1)-(3) in Proposition \ref{added1}.  A totally different
approach to guaranteeing this convergence, which  allows general
continuous $\ell$,  is  as follows. Set
\begin{equation}
\label{infff} V^r_\infty(x):= \inf_{\alpha\in {\cal A}^r}
J[\ell,h](x,+\infty,\alpha)\; \; \in\; \; [-\infty,+\infty],
\end{equation}
with the same conventions (i)-(ii) used to define $V_\infty$ in
(\ref{inff}).  \bp{added2} Assume $(A_1)$-$(A_5)$ hold  and $\cg$
is asymptotically null for $f^r$. Let $w:\cg\to\R$ be a continuous
solution of (\ref{HJBE}) on $\cg\setminus\{0\}$ satisfying
$w(0)=0$.  Then $w\equiv V^r_\infty$ on $\cg$. \ep \bpr We show
how to modify the proof of Theorem \ref{thm1}.  We again assume
$h\equiv 0$, the general case being handled in a similar way. Let
$\bar x\in \cg\setminus\{0\}$. For each $\eps>0$, there exists
$(1,\alpha_\eps)\in {\cal Z}_1$ such that
\begin{eqnarray}\label{ipps} w(\bar x) &\ge&
\int_0^1\ell(\phi(s,\bar x,\alpha_\eps),\alpha_\eps(s))\, {\rm
d}s\\&+&w(\phi(1,\bar x,\alpha_\eps))-\eps \ge  w(\bar x)-
\eps\nonumber\end{eqnarray} with the last inequality following
from the first part of Lemma \ref{lemma2}. Using Lemma
\ref{compactness} with $c=1$, we can find a subsequence of the
$\alpha_\eps$'s (which we do not relabel) and $\beta\in {\cal
A}^r$ such that (i) $\alpha_\eps \to\beta$  weak-$\star$ on
$[0,1]$ and (ii) $\max\{||\phi(t,\bar x,\alpha_\eps)-\phi(t,\bar
x,\beta)||: 0\le t\le 1\}\to 0$ as $\eps\downarrow 0$. Since $\cg$
is $f^r$-invariant, it follows that $\phi(t,\bar x,\beta)\in \cg$
for all $t\in [0,1]$.  Letting $\eps\downarrow 0$ in (\ref{ipps})
and using the continuity of $w$ now gives \[
 w(\bar x)=
\int_0^1\ell^r(\phi(s,\bar x,\beta),\beta(s))\, {\rm
d}s+w(\phi(1,\bar x,\beta))
\]
The proof of the preceding equality is based on the continuity of
the maps $a\mapsto \ell(x,a)$ for all $x$, and is similar to the
argument in the appendix of \cite{M00a}.  This procedure is
iterated and gives $\hat\alpha\in {\cal A}^r$ such that
\begin{equation}\label{esa}
w(\bar x)= \int_0^M\ell^r(\phi(s,\bar
x,\hat\alpha),\hat\alpha(s))\, {\rm d}s+w(\phi(M,\bar
x,\hat\alpha))\; \; \; \; \forall M\in\N
\end{equation}
Since $\cg$ is asymptotically null for $f^r$ and $w$ is
continuous, $w(\phi(M,\bar x,\hat\alpha))\to 0$ as $M\to +\infty$.
This implies that the integral in (\ref{esa}) converges in $\R$ as
$M\to+\infty$, so (\ref{esa}) gives $w(\bar x) \ge V^r_\infty(\bar
x)$.  The proof that $w(\bar x)\le V^r_\infty(\bar x)$ is similar
to  the proof of the corresponding inequality in the proof of
Theorem \ref{thm1}, with ${\cal A}$ replaced by ${\cal A}^r$,
since (cf. \cite{BCD97})
\[\sup\{-f^r(x,a)\cdot p-\ell^r(x,a): a\in A^r\}\; =\;
\sup\{-f(x,a)\cdot p-\ell(x,a): a\in A\}\] for all $x\in\R^N$.
 \epr
\br{added3} In some applications, we will take $f$  ULES, $h\equiv
0$, and $\ell=-g$ for uniformly locally Lipschitz $g$ (cf. $(A_6)$
below). In this case, $V_\infty$ is a continuous solution of
(\ref{HJBE}) on ${\cal D}^r_o\setminus\{0\}$ (with $h\equiv 0$),
where ${\cal D}^r_o$ is the robust domain of attraction for
$f^r:\R^N\times A^r\to\R^N$.\footnote{Recall that the {\em relaxed
domains of attraction} are defined in the following way: ${\cal
D}^r:=\left\{x\in\R^N: \ds\lim_{t\to+\infty}\phi(t,x,\alpha)= 0 \;
\forall \alpha\in {\cal A}^r\right\}$, and $\ds{\cal
D}^r_o:=\left\{x\in \R^N: \sup_{\alpha\in {\cal
A}^r}t(x,\alpha)<+\infty\right\}$.} Proposition \ref{added2} then
implies that $V_\infty\equiv V^r_\infty$ on ${\cal D}^r_o$.
 We will show
below that we in fact have $V_\infty\equiv V^r_\infty$ on all of
${\cal D}_o$, since ${\cal D}_o={\cal D}^r_o$.
 \er
\section{Properties of Domains of Attraction and \- Lyapunov Functions} \label{appl}
Our cost function $g$ gives rise to value functions, which we
denote by $V_L$, for infinite horizon cost maximization (cf.
(\ref{crit})).  In this section, we add hypotheses on $g$ which
imply that $V_L$ is a robust Lyapunov function for $f$, and that
$V_L$ is a unique solution of the PDE (\ref{auxx}).  We also give
PDE characterizations for ${\cal D}_o$.
 Regularity for
$V_L$ is covered in $\S$\ref{appl2}. In particular, we show that
by removing the positive lower bound assumption (\ref{poss}) on
$g$, we obtain Lyapunov functions with special properties   not
found in the Lyapunov functions of \cite{CGW00,CGW00a, CGW01}. In
$\S$\ref{appl3}, these results are used to give uniqueness
characterizations for solutions of the generalized Zubov equation
(\ref{gz}), giving  a new, more general class of robust
 Lyapunov functions which can be
characterized as unique solutions of generalized Zubov equations.
\subsection{Domains of Attraction}
To study ${\cal D}_o$, we first consider the {\em auxiliary} PDE
\begin{equation}\label{negg}
\ds \inf_{a\in A}\{-f(x,a)\cdot Dv(x)-g(x,a)\}=0
\end{equation}
under hypotheses $(A_1)$-$(A_4)$.  It is immediate from the
definition of viscosity solutions that a function $w$ is a
viscosity solution of (\ref{negg}) on an open set ${\cal O}$
exactly when $-w$ is a viscosity solution of (\ref{HJBE}) on
${\cal O}$, with the choices $\ell=-g$ and $h\equiv 0$.  Equation
(\ref{HJBE}) is  the usual Hamilton-Jacobi equation, which was
studied for general nonnegative $\ell$ in  \cite{M00a, M00}.
However, it will be convenient for our applications to assume that
$g$ is nonnegative, so the earlier uniqueness results on
(\ref{HJBE}) for nonnegative $\ell$ and $h\equiv 0$ would  not
apply. Instead, we study (\ref{negg}) using the theory of the
previous section. We set
\begin{equation}
\label{crit} V_L(x)=\sup_{a\in {\cal A}} J[g](x,+\infty,a)\in
[0,+\infty],\; \; x\in\R^N
\end{equation} We sometimes write $V_L[g]$ instead of
$V_L$, to emphasize the cost function $g$. If $(A_1)$-$(A_4)$ hold
with $f$ ULES, and if
\begin{equation}\label{rella}
  \exists \; \; \tilde C,\; \lambda>0 \; \; \text{s.t.} \; \;
  g(x,a)\le \tilde
  C||x||^\lambda \; \; \text{for\ \  all}\; \;
   x\in B_r \; \; \text{and} \; \; a\in
  A
\end{equation}
then the argument of \cite{CGW01} shows that $V_L[g]$ is locally
bounded on ${\cal D}_o$ (i.e., $\sup\{V_L(x): x\in K\}<\infty$ for
each compact set $K\subseteq {\cal D}_o$).  Moreover, if $g$
satisfies the {\em stronger condition}\bi\item[]\bi
\item[$(A_6)$\ \ ]For each $R>0$,  there exists a positive
constant  $L_{g,R}$ such that $|g(x,a)-g(y,a)|\le L_{g,R}||x-y||$
for all $x,y\in B_R$ and $a\in A$. \ei\ei then $V_L[g]$ is
continuous on ${\cal D}_o$ (cf. \cite{CGW01}).
 The
following theorem extends the uniqueness result in \cite{CGW01},
Theorem 3.9, to the case of general nonnegative cost functions:
\bt{thm2} Assume $(A_1)$-$(A_4)$ and (\ref{rella}) hold and $f$ is
ULES. Let ${\cal O}\subseteq\R^N$ be a $g$ quasi-stable set
containing the origin, and let $w:{\cal O}\to\R$ be a
bounded-from-below solution of (\ref{negg}) on ${\cal O}$
satisfying condition $(SC_w)$ and $w(x)\to +\infty$ as $x\to x_o$
for all $x_o\in\partial{\cal O}$. Then ${\cal O}={\cal D}_o$, and
$w\equiv V_L$ on ${\cal D}_o$.\et\bpr We can assume that $w$ is
continuous on $\calo$. (Indeed, in what follows, we do not use the
full strength of the continuity of $w$.  Instead, we only use the
lower semicontinuity of $w_\star$ and the upper semicontinuity of
$w^\star$.) As already noted, $-w$ is a solution of (\ref{HJBE})
with $\ell=-g$ and $h\equiv 0$ on ${\cal O}$. It follows from
Lemma \ref{lemma2} and the boundary behavior of $w$ that ${\cal
O}$ is $f$-invariant. Indeed, if $x\in {\cal O}$ and $\alpha\in
{\cal A}$, and if $t$ is the first time $\phi(t,x,\alpha)\in
\partial{\cal O}$, then we can apply the first part of
Lemma \ref{lemma2} to (\ref{HJBE}), with $q=x$,  $\ell=-g$,
$h\equiv 0$, and $\beta=\alpha$, and with bounded open sets
$B\subseteq \calo$ containing ${\rm
Trace}(\phi(\cdot,x,\alpha)\lceil[0,t-1/m])$ for $m\in \N$,  to
get
\[\begin{array}{l}
-w(x)\le -\int_0^{t-1/m}g(\phi(s,x,\alpha),\alpha(s))\, {\rm
d}s-w(\phi(t-1/m,x,\alpha))\to -\infty\end{array}
\]
as $m\to +\infty$,   which is a contradiction.  Also, the first
part of Lemma \ref{lemma2} and the invariance of $\calo$ then give
\begin{equation}
\label{gives}
 w(x)\ge \int_0^tg(\phi(s,x,\alpha),\alpha(s))\, {\rm
d}s +w(\phi(t,x,\alpha))
\end{equation}
for all $\alpha\in {\cal A}$, $x\in \calo$, and $t>0$. Since $w$
is bounded-from-below and $g$ is nonnegative, (\ref{gives}) gives
$J[g](x, +\infty,\alpha)<\infty$ for all $\alpha\in {\cal A}$ and
$x\in \calo$. The $g$ quasi-stability of ${\cal O}$ and Lemma
\ref{lemma1} therefore give ${\cal O} \subseteq {\cal D} \subseteq
{\rm cl}({\cal D}_o)$. Since $\calo$ is open, $\calo\subseteq
{\cal D}_o$. Since ${\cal D}_o$ is asymptotically null for $f$ and
$\calo$ is an open invariant set containing the origin, $\calo$ is
also asymptotically null for $f$. Therefore,  Theorem \ref{thm1}
(with the choices $\ell=-g$ and $h\equiv 0$) gives $w\equiv V_L$
on $\calo$. Therefore, $V_L(x)\to+\infty$ as $\calo\ni x\to x_o$
for all $x_o\in
\partial \calo$.  By (\ref{rella}),
$V_L$ is locally bounded on ${\cal D}_o$, so we conclude that $
{\cal D}_o\cap
\partial\calo=\emptyset
$.
 Recall that $\calo$ contains a
neighborhood of the origin. If $p_o\in {\cal D}_o\setminus\calo$,
then the asymptotic nullness of ${\cal D}_o$ gives a point $p_1\in
{\cal D}_o\cap
\partial\calo$, namely, the first point on
a trajectory for $f$ starting at $p_o$ and approaching the origin
that lies in $\partial\calo$. This contradiction gives ${\cal
D}_o\subseteq \calo$ and completes the proof. \epr

Notice that {\em any} open set $\calo\subseteq\R^N$ is $g$
quasi-stable if (i) $f$ is ULES and (ii) $g$ satisfies
(\ref{poss}) for the constant $r$ in the ULES definition, which is
the lower bound requirement  in \cite{CGW00,CGW00a, CGW01}.
Therefore, the uniqueness results for (\ref{negg}) given in
\cite{CGW01} all follow from Theorem \ref{thm2}.  Moreover, the
comparison results for (\ref{negg}) given in \cite{CGW01} also
require the comparison functions $w$ to be proper (meaning
$w(x)\to +\infty$ as $||x||\to\infty$) and $(A_6)$.  Therefore,
Theorem \ref{thm2} {\em extends} the uniqueness results in
\cite{CGW01} by taking a more general class of comparison
functions into account. Moreover, Theorem \ref{thm2} applies to
more general $g$, including non-Lipschitz $g$ that do not satisfy
(\ref{poss}) (cf. $\S$\ref{illu}). If $(A_1)$-$(A_4)$ and $(A_6)$
hold, with $f$ ULES,  then $V_L$ is a continuous solution of
(\ref{negg}) on ${\cal D}_o$ by arguments in \cite{CGW01}, so
Theorem \ref{thm2} gives PDE characterizations for $V_L$ and
${\cal D}_o$ if $V_L(x)\to +\infty$ as $x\to\partial ({\cal
D}_o)$. Sufficient conditions for this boundary behavior of $V_L$
will be given in the next subsection.
\subsection{Further Properties of $V_L$}
\label{appl2} In this section, we study the behavior of $V_L$ near
$\partial({\cal D}_o)$ and the regularity of $V_L$.   If we assume
$(A_1)$-$(A_4)$ and $(A_6)$ with $f$ ULES, and if we also
assume\bi
\item[$(A_7)$]\ \
$\int_0^tg(\phi(s,x,a),a(s))\, {\rm d}s>0$ for all $x\in {\cal
D}_o \setminus\{0\}$, $t>0$, and $a\in {\cal A}$ \ei then $V_L$ is
a robust Lyapunov function for $f$ on ${\cal D}_o$. Indeed, the
dynamic programming methods from \cite{BCD97} give
\begin{equation}
\label{dpp} \ds V_L[g](x)=\sup_{\alpha\in {\cal A}} \left[
\int_0^tg(\phi(s,x,\alpha),\alpha(s))\, {\rm
d}s+V_L[g](\phi(t,x,\alpha))\right]\; \; \forall t>0
\end{equation}
for all $x\in{\cal D}_o$, which gives
$V_L(x)>V_L(\phi(t,x,\alpha))$ for all $x\in {\cal
D}_o\setminus\{0\}$, $t>0$, and $\alpha\in {\cal A}$, as needed.
Since Lyapunov functions form the basis for much of current work
in stability, this motivates our study of further properties of
$V_L$. Notice that $(A_7)$ is weaker than (\ref{poss}) (cf. Remark
 \ref{compare} and  $\S$\ref{illu} below for examples).
\subsubsection{Boundary Behavior of $V_L$}
We will assume $(A_1)$-$(A_4)$ and $(A_6)$-$(A_7)$, with $f$ ULES,
throughout the remainder of this subsection. We first look for
conditions under which
\begin{equation}\label{asy}
  V_L(x)\; \to\;  +\infty\; \; \text{as}\; \; x\to x_o\; \;
  \text{for\ \  all}\; \; x_o\in \partial({\cal D}_o)
\end{equation}
Condition (\ref{asy}) clearly holds if there is a function $w$
satisfying the hypotheses of Theorem \ref{thm2}. In \cite{CGW01},
$g$ is assumed to satisfy (\ref{poss}), which implies (\ref{asy}),
because $V_L(x)\ge g_o[\sup_\alpha t(x,\alpha)]\to +\infty$ as
$x\to \partial ({\cal D}_o)$. More generally, (\ref{asy}) will
hold if we add the assumption \bi
\item[$(A_8)$]\ \
$\R^N$ is  $g$ quasi-stable \ei We say that $g$ satisfies the {\bf
standing hypotheses} if $(A_1)$-$(A_4)$ and $(A_6)$-$(A_8)$ hold.
We now show that (\ref{asy}) holds under the added hypothesis
$(A_8)$. Set
\[
V^r_L(x)=\sup_{a\in {\cal A}^r}J[g](x,+\infty,\alpha)\; \in\; [0,
+\infty]
\]
 We sometimes
write $V^r_L[g]$ instead of  $V^r_L$, to emphasize the cost
function $g$. Recall (cf. \cite{BCD97}) that (\ref{negg}) has
exactly the same set of solutions on any open set as
\begin{equation}
\label{rbell} \ds\inf_{a\in A^r}\{-f^r(x,a)\cdot
Dv(x)-g^r(x,a)\}=0,
\end{equation}
and that $V_L^r$ is a continuous  solution of (\ref{rbell}) on
${\cal D}^r_o$, by \cite{CGW01}.   Recall that $f^r$ is ULES.  By
Lemma \ref{lemma1}, we know that ${\cal D}^r_o$ is open  and that
${\cal D}_o^r={\cal D}^r$, by the convexity of the sets
$f^r(x,A^r)$. Notice that $(A_8)$ and the equivalence of $g$
quasi-stability and relaxed $g$ quasi-stability  (cf. Appendix A)
imply that ${\rm dom}(V_L^r)\subseteq {\cal D}^r$. Therefore,
${\rm dom}(V_L^r) = {\cal D}^r_o$, since $V_L^r$ is finite on
${\cal D}^r_o$ (cf. \cite{CGW01}). By Theorem \ref{thm2},
(\ref{asy}) will follow if
\begin{equation}
\label{relaxed} V^r_L(x)\to +\infty \; \; \text{as}\; \; x\to x_o
\; \; \text{for \ all}\; \; x_o\in\partial \left({\cal
D}_o^r\right),
\end{equation}
 since then we can take $w=V^r_L$ and $\calo={\cal D}_o^r$
 in Theorem \ref{thm2}
and conclude that ${\cal D}_o^r={\cal D}_o$ and $V_L^r=V_L$ on
${\cal D}_o$.
 We will now prove (\ref{relaxed}).

 Suppose
  that $x_n\to x_o\in\partial({\cal D}^r_o)$, but that there
  exists $K$ such that
$V^r_L(x_n)\le K<\infty$ for all $n$.  Fix $a\in {\cal A}^r$ and
$M>0$.  Set \[{\cal S}_M=B_1(\{ \phi(t,x_o,a): 0\le t\le M\}),\]
and choose $\bar F>0$ such that $||f(x,a)-f(y,a)||\le \bar F
||x-y||$ for all $x$ and $y$ in $B_1({\cal S}_M)$ and all $a\in
A$.
 There exists an $N_M\in \N$ such
that
\begin{equation}\label{gron1}
||\phi(t,x_o,a)-\phi(t,x_n,a)||\; \; \le\; \; e^{t {\bar F}
}||x_n-x_o||\; \; \; \; \forall t\in [0,M],\; \; \forall n\ge N_M,
\end{equation}
and such that $\phi(t,x_n,a)\in {\cal S}_M$ for all $t\in [0,M]$
and all $n\ge N_M$.    This follows from Gronwall's Inequality and
a generalization of standard estimates from \cite{BCD97}.  By
enlarging $N_M$, we can guarantee that
\begin{equation}\label{gron2}
\ds ||x_n-x_o||\; \; \le\; \; \frac{1}{M\left[{\bar g}
+1\right]}e^{-M {\bar F}}\; \; \; \; \forall n\ge N_M
\end{equation}
where  $|g(x,a)-g(y,a)|\le \bar g||x-y||$ for all $x,y\in
B_1({\cal S}_M)$ and $a\in A$.  Combining (\ref{gron1}) and
(\ref{gron2}),
\begin{equation}
\label{integrand}
 g^r(\phi(t,x_o,a),a(t))\; \; \le\; \;
g^r(\phi(t,x_n,a), a(t))\; +\; 1/M\; \; \; \; \forall n\ge N_M
\end{equation}
a.e. $t\in[0,M]$. We can integrate (\ref{integrand}) to get
\begin{equation}\label{inte}
J[g](x_o,M,a)\le K+1\; \; \forall  a\in {\cal A}^r,M>0.
\end{equation} Letting
$M\to +\infty$ in (\ref{inte}) for fixed $a$, and using the
equivalence of $g$ quasi-stability and relaxed $g$
quasi-stability, we conclude that $x_o\in {\cal D}^r$. Recalling
that $f^r(x,A^r)$ is convex for all $x$, Lemma \ref{lemma1}(c)
gives $x_o\in {\cal D}_o^r$, which contradicts the fact that
${\cal D}_o^r$ is open. The following proposition summarizes these
observations: \bc{prop1} Let $(A_1)$-$(A_4)$ and $(A_6)$-$(A_8)$
hold with $f$ ULES.  Then $V_L$ is a robust Lyapunov function for
$f$ on ${\cal D}_o$ which satisfies (\ref{asy}).   Moreover,
${\cal D}_o={\cal D}^r_o$, and $V_L\equiv V^r_L$ on ${\cal D}_o$.
\ec

\br{rki} Under our standing hypotheses,  \[{\cal D}_o\subseteq
{\rm dom} (V_L)\subseteq {\cal D} \subseteq {\rm cl}({\cal
D}_o),\] which is weaker than the condition ${\rm dom}(V_L)={\cal
D}_o$. Later, we show that under $(A_1)$-$(A_8)$, we do in fact
have ${\rm dom} (V_L)={\cal D}_o$. Our hypotheses are {\em weaker}
than (\ref{poss}), since $(A_7)$ allows vanishing  and
asymptotically decaying costs $g$ (cf. \cite{M00},  and Remark
\ref{compare} and $\S$\ref{illu} below). For the study of the case
where (\ref{poss}) holds, see \cite{CGW00, CGW00a, CGW01, GW00}.
Notice that $(A_1)$-$(A_8)$ imply that $\R^N$ is also relaxed $g$
quasi-stable (cf. Appendix A). Also, the proof of  (\ref{asy})
remains valid if $(A_8)$ is relaxed to the requirement that there
be a $g$ quasi-stable set containing ${\rm cl}({\cal D}_o)$.\er

\br{prop} Under the additional hypothesis (\ref{poss}), $V_L[g]$
is proper (cf. \cite{CGW01}).  However, under the standing
hypotheses on $g$, $V_L[g]$ may not be proper.  For example, we
can take \[N=1,\; \;  g(x)=|x|/(1+x^2),\; \;  \text{and}\; \;
f(x,a)\equiv -x,\] in which case $V_L[g](x)\to \pi/2$ as
$x\to+\infty$.\er

\subsubsection{Lipschitz Lyapunov Functions}
\label{lipsec} It is natural to look for conditions under which
the Lyapunov function $V_L$ is not only continuous on ${\cal
D}_o$, but also locally Lipschitz on ${\cal D}_o$.  It is also
natural to ask whether the function
\[
v_\delta:=1-e^{-\delta V_L[g]}
\]
is globally Lipschitz for large enough fixed $\delta>0$, because
$v_\delta(x)\to 1$ as $\delta\to +\infty$
 for each fixed $x\in
{\cal D}_o\setminus\{0\}$.
  Later, we will show that $v_\delta$ is a robust Lyapunov
function for $f$ on ${\cal D}_o$ for each fixed $\delta>0$, and
that it is a global solution of the Zubov equation for suitable
$g$.  The question of whether $v_\delta$ is globally Lipschitz for
fixed $\delta>0$ therefore reduces to the question of whether the
methods of the previous section can produce globally Lipschitz
Lyapunov functions and globally Lipschitz solutions of (\ref{gz}).
The proof of the following is the same as the proofs of
Proposition 4.2 and Proposition 4.3 in \cite{CGW01}: \bp{lift2}
Let $(A_1)$-$(A_4)$ and $(A_6)$-$(A_7)$ hold with $f$ ULES, and
assume the following: \bi
\item[(L)]
There are   $c,\kappa>0$ and $s>L_r/\sigma$ such that for all
$x,y\in B_c$ and $a\in A$, \[|g(x,a)-g(y,a)|\le \kappa
\max\{||x||^s,||y||^s\}||x-y||,\] where $r$ and $\sigma$ are the
constants in the ULES definition and $L_r$ is the constant  from
$(A_2)$. \ei
 Then
$V_L$ is locally Lipschitz on ${\cal D}_o$.  If in addition  there
are $L_f\in (0,s\sigma)$ and $L_g>0$ such that
\begin{equation}
\label{globallip} ||f(x,a)-f(y,a)|| \le L_f||x-y|| \; \;
\text{and}\; \; |g(x,a)-g(y,a)| \le L_g||x-y||
\end{equation} for all $x,y\in\R^N$ and $a\in A$,  and if (\ref{poss}) holds,
then $v_\delta$ is globally Lipschitz on $\R^N$ for all
sufficiently large $\delta>0$. \ep

The proof of the global Lipschitzness of $v_\delta$ is based on
the fact that $V_L(x)>g_o [\sup\{t(x,a): a\in {\cal A}\}]$, where
$g_o$ is the uniform lower bound for $g$ in (\ref{poss}).  The
following example shows that (\ref{poss}) cannot be deleted from
the hypotheses of this proposition.  It illustrates how relaxing
the uniform lower bound on $g$ gives  phenomena not found in the
Zubov equation solutions in \cite{CGW01} (see also Remark
\ref{rknewly}).

\bex{hav} For each $k\in \N_o:=\{0,1,2,\ldots\}$, set
\[
t_{k-}=10^k-\frac{1}{10^{2k+1}}\; \; \text{and}\; \;
t_{k+}=10^k+\frac{1}{10^{2k+1}}.
\]
Define $I_\Delta$ and $\Delta:I_\Delta\to\R$ by
\[\ds
I_\Delta=\bigcup_{k\in \N_o}\left[t_{k-}, t_{k+}\right],\;
\Delta(x)= \twoif{10^{3k+1}\left(x-t_{k-}\right),}{t_{k-}\le x\le
10^k, \; k\in \N_o}{10^{3k+1}\left(t_{k+}-x\right),}{10^k\le x\le
t_{k+}, \; k\in \N_o}\] Then the graph of $\Delta$ is a sequence
of nonoverlapping triangles centered at the points $10^k$ for
$k\in \N_o$ which become taller and thinner as $k\to +\infty$. In
fact, while \[\Delta(10^k)=10^k\] for all $k\in \N_o$, we have
\[
\int_{I_\Delta} \Delta(x)\, {\rm d}x=\frac{1}{10}\sum_{k=0}^\infty
10^{-k}<\infty
\]
Let $Q:\R\to \R$ be any continuous function that satisfies the
following conditions. \bi
\item
$Q\equiv \Delta$ on $I_\Delta$, $Q(x)= (\frac{9}{10})^6x^6$ for
$0\le x\le \frac{9}{20}$, and $Q(x)\ge 0$ for all $x\ge 0$
\item
${\rm Null}(Q)=\{0, \pm t_{1-}, \pm t_{1+}, \pm t_{2-}, \pm
t_{2+}, \ldots\}$, and $|Q(x)|\le 1$ for all $x$ in ${\rm
comp}(I_\Delta)$
\item
$Q\lceil[{\rm comp}\left(I_\Delta\right)]$ is Lipschitz with
Lipschitz constant ${\cal L}\le 1$,  $Q$ is odd, and
$\int_{[0,\infty)} Q(x)\, {\rm d}x<\infty$ \ei We leave the easy
construction of $Q$ to the reader. The functions
\[
g(x)=-Q(x)f(x),\;  \; \; \; \; \;
f(x)=-\ds\twoif{\left(\frac{10}{9}\right)^6x,}{-\frac{9}{10}\le
x\le \frac{9}{10}}{\frac{1}{x^5},}{x\ge \frac{9}{10}\; \text{or}\;
x\le -\frac{9}{10}}
\]
satisfy (\ref{globallip}) with $L_g>0$ and $L_f\in
(0,6(\frac{10}{9})^6)$. Take $f$ as the dynamics, with no
 controls, and $g$ as the cost function. Then
$(A_1)$-$(A_4)$ and $(A_6)$-$(A_7)$ hold with $f$ ULES, and $(L)$
holds with parameters $ r=9/10$, $C=1$,
$\sigma=L_r=\left(10/9\right)^6$, $s=6$, $\kappa=7$, and $c=9/20$.
We let $\phi(t,x)$ denote the trajectory for $f$ and the initial
value $x$.  For all $x>0$, it follows that
\[
V_L(x)=\int_0^\infty g(\phi(t,x))\, {\rm d}t=\int_0^\infty
\frac{g(\phi(t,x))}{f(\phi(t,x))}\, \frac{\partial\phi}{\partial t
}(t,x)\, {\rm d}t=\int_0^xQ(u)\, {\rm d}u.\] Therefore, if
$\delta>0$ is given, and if we set $v_\delta(x)=1-e^{-\delta
V_L(x)}$, then we get
\[
|Dv_\delta(10^k)|=\delta Q(10^k)\exp\left(-\delta
\int_0^{10^k}Q(s)\, {\rm d}s\right)\to +\infty\; \; \text{as}\; \;
k\to +\infty.\] Therefore, while $V_L$ is locally Lipschitz, there
cannot exist $\delta>0$ such that $v_\delta$ is  globally
Lipschitz.  This example shows that the positivity condition
(\ref{poss}) cannot be omitted from the statement of Corollary
\ref{lift2}.
 \eex

\subsection{Uniqueness of Viscosity Solutions of Generalized Zubov Equation}
\label{appl3} For any function $w:\R^N\to(-\infty, +\infty]$, we
define the  {\bf Kru\v{z}kov transformation} $\kw:\R^N\to \R$ of
$w$ by $\check{w}(x):=1-e^{-w(x)}$, with the convention that
$e^{-\infty}=0$.  If $w$ is a solution of the auxiliary PDE
(\ref{negg}) on ${\cal D}_o$, then  $\check w$ is a solution of
the generalized Zubov equation (\ref{gz}) on ${\cal D}_o$ (cf.
\cite{BCD97}, Chapter II). We apply this observation to $V_L$,
assuming for the rest of this subsection that the hypotheses of
Corollary \ref{prop1} are satisfied.

Recall that ${\rm dom}(V^r_L)= {\cal D}_o^r(={\cal D}^r)$. Since
$V^r_L(x)\to +\infty$ as $x\to\partial({\cal D}^r_o)$  and $V^r_L$
is continuous on ${\cal D}^r_o$,  $\check V^r_L$ is continuous on
$\R^N$. One can check (cf. \cite{CGW01}) that $\check V^r_L$
satisfies the Dynamic Programming Principle
\begin{equation}
\label{DP}
\begin{array}{l}
\ds \check V^r_L(x)=\sup_{\alpha\in {\cal A}^r}\left\{
[1-G(x,t,\alpha)]+G(x,t,\alpha) \check
V^r_L(\phi(t,x,\alpha))\right\}\; \; \forall t\ge 0,\\
G(x,t,\alpha):=e^{-J[g](x,t,\alpha)}\end{array}
\end{equation} on $\R^N$.  Using (\ref{DP}),
one can use the arguments of \cite{BCD97} to show that $\check
V^r_L$ is a viscosity solution of the generalized Zubov equation
(\ref{gz}) on all of $\R^N$. Also, (\ref{DP}) and $(A_7)$ give
\begin{eqnarray}
\check V_L(\phi(t,x,\alpha))&\le&
\left[1-G(x,t,\alpha)\right]\left[\check
V_L(\phi(t,x,\alpha))-1\right]+\check V_L(x)\; <\;  \check
V_L(x)\nonumber\end{eqnarray} for all $x\in {\cal
D}_o\setminus\{0\}$, $\alpha\in {\cal A}$ and $t>0$,
 since
$V^r_L=V_L$ on ${\cal D}_o$ (cf. $\S$\ref{appl2}).   Therefore,
$\check V_L$ is also a robust Lyapunov function for $f$ on ${\cal
D}_o$.

Now let $w:\R^N\to\R$ be any solution of (\ref{gz}) on $\R^N$ that
satisfies $(SC_w)$.  By applying Theorem \ref{thm1} with
$\ell\equiv -g$ and $h\equiv g$, we conclude that $w$ agrees with
\begin{equation}
\label{dough} \check V_L(x)\; \; =\; \; \sup_{\alpha\in {\cal A}}
\int_0^\infty G(x,t,\alpha)\,  g(\phi(t,x,\alpha),\alpha(t))\,
{\rm d}t\; \; \in \R
\end{equation}
on the asymptotically null set ${\cal D}_o$. In fact, if we also
assume that $w$ is bounded, then $w\equiv \check{V}_L$ on all of
$\R^N$. To see why, we can assume $w$ is continuous, the general
case being proven in a similar way. First note that by the
boundedness of $w$ and $(SC_w)$,
\begin{equation}\label{lom}
 \lim_{t\to +\infty}G(x,t,\alpha)\, w(\phi(t,x,\alpha))=0\; \; \; \;
\forall \alpha\in {\cal A},\; \forall x\in\R^N
\end{equation}
which follows by separately considering the cases where the
exponent  \[\int_0^\infty g(\phi(s,x,\alpha),\alpha(s))\, {\rm
d}s\] converges or diverges and using  $g$ quasi-stability. An
application of the first part of Lemma \ref{lemma2}, $(SC_w)$, and
(\ref{lom}) then gives $w(x)\ge J[g,g](x,+\infty,\alpha)$ for all
$\alpha\in {\cal A}$ and $x\in \R^N$, so $w(x)\ge \check V_L(x)$.
On the other hand, given $\eps>0$ and $x\in\R^N$, the construction
in the proof of Theorem \ref{thm1}  gives an input $\hat\alpha\in
{\cal A}$ such that
\begin{eqnarray}
-w(x) &\ge& -\int_0^tG(x,s,\hat\alpha)\,
g(\phi(s,x,\hat\alpha),\hat\alpha(s))\, {\rm d}s\nonumber\\&-&
G(x,t,\hat\alpha)\,
w(\phi(t,x,\hat\alpha))-\eps(1-e^{-t})\nonumber\end{eqnarray} for
all $t\in \N$. Combining this and (\ref{lom}) now gives \[w(x)\;
\le\; J[g,g](x,+\infty,\hat\alpha)+\eps\; \le\;  \check
V_L(x)+\eps,\] as needed. Recall that $\check V^r_L$ is a
continuous bounded viscosity solution of (\ref{gz}) on $\R^N$
which is null at the origin. Setting $w=\check V^r_L$ in the
preceding argument therefore allows us to conclude that $\check
V^r_L=\check V_L$ on all of $\R^N$, so $\check V_L$ is also a
continuous solution on $\R^N$. The following corollary summarizes
these observations and includes the global uniqueness results for
the generalized Zubov equation in \cite{CGW01}. \bc{prop2} Assume
$(A_1)$-$(A_4)$ and $(A_6)$-$(A_8)$ with $f$ ULES. Then:
\bi\item[(1)] $\check V_L$ is a robust Lyapunov function for $f$
on ${\cal D}_o$, and  ${\cal D}_o=\check V_L^{-1}([0,1))$.
\item[(2)]
If $w$ is a solution of (\ref{gz}) on ${\cal D}_o$ satisfying
$(SC_w)$, then $w\equiv \check V_L$ on ${\cal D}_o$.
\item[(3)]
If $w$ is a bounded solution of (\ref{gz}) on $\R^N$ satisfying
$(SC_w)$, then $w\equiv \check V_L$ on $\R^N$. \ei Moreover,  (i)
$\check V_L$ is the unique bounded solution $w$ of (\ref{gz}) on
$\R^N$ that satisfies $(SC_w)$ and (ii) ${\cal D}_o=w^{-1}([0,1))$
for any bounded solution $w$ of (\ref{gz}) on $\R^N$ that
satisfies $(SC_w)$.\ec

\br{recall} In Example \ref{ex1}, $(A_1)$-$(A_4)$ and
$(A_6)$-$(A_7)$ hold, but $(A_8)$ is not satisfied, and there are
{\em infinitely many} bounded solutions of (\ref{gz}) on $\R^N$
satisfying $(SC_w)$. Therefore, the quasi-stability hypothesis
$(A_8)$ of Corollary \ref{prop2} cannot be dropped.\er

\br{rknewly} Conclusions (2)-(3) of the preceding corollary remain
true if $(A_6)$ is replaced by (\ref{rella}).  We can also prove
`nonglobal' PDE characterizations for solutions of (\ref{gz}) on
general open sets $\calo$ under our relaxed conditions on $g$ (cf.
Appendix B).  The paper \cite{GW00} suggests the problem of
determining what subset of the set of all robust Lyapunov
functions $V$ for a given ULES dynamics $f$ has the following
properties: (i) $V\equiv \check V_L[g]$ for some cost function $g$
and (ii) $V$ is the
 unique bounded solution of the generalized Zubov
equation (\ref{gz}) that is null at the origin. One consequence of
our results is that by allowing more general  costs $g$, we made
the subset of known Lyapunov functions satisfying these two
properties strictly larger. In particular, the set of Lyapunov
functions $\check V_L[g]$ studied in \cite{CGW01, GW00} is a {\em
proper} subset of the set of all functions that can be written as
unique solutions of generalized Zubov equations. To see why, we
have to find robust Lyapunov functions $V$ which can be written as
$\check V_L[g]$ for cost functions $g$ satisfying our standing
hypotheses, but which cannot be written as $\check V_L[\hat g]$ if
$\hat g$ is also assumed to satisfy the positivity condition
(\ref{poss}) which is assumed in \cite{CGW00, CGW00a, CGW01,
GW00}. A general method for doing this is as follows.

Let $f$ be a ULES dynamics satisfying $(A_1)$-$(A_3)$ and
$B_r\subsetneq {\cal D}_o$, and $g$ satisfy the standing
hypotheses  in such a way that $g(x,a)\equiv \gamma(x)$ and
$\gamma(\bar x)=0$ for some $\bar x\in {\cal D}_o\setminus B_r$
(cf. \S\ref{illu} for particular cases).  In particular, this
means $\gamma$ is locally Lipschitz (cf. $(A_6)$). Then $g$ does
not satisfy the positivity condition (\ref{poss}). Take $W$ to be
the associated Lyapunov function $V_L[g]$, so that
$1-e^{-W}=\check V_L[g]$. Then $W$ is continuous on ${\cal D}_o$,
and $\cvlg$ is the unique bounded solution of (\ref{gz}) on $\R^N$
satisfying $(SC_w)$. Suppose $\hat g$ satisfies the standing
hypotheses and also the positivity condition (\ref{poss}), and
suppose further  that $W\equiv V_L[\hat g]$ on ${\cal D}_o$. Let
$L_\gamma>1$ be a Lipschitz constant for $\gamma$ on $B_1(\bar
x)$, and pick $\eps\in (0,g_o)$, where $g_o$ is   from the
positivity  condition (\ref{poss}) on $\hat g$.  Since $f$ is
bounded, we can find $t>0$ such that
\[
||\phi(s,\bar x,\alpha)-\bar x||\le \frac{\eps}{L_\gamma}\wedge
\left[ \frac{1}{2} \dist(\bar x,B_r\cup \partial ({\cal
D}_o))\right]\wedge 1/2\; \; \; \; \; \forall s\in [0,t],\;
\forall \alpha\in {\cal A}^r
\]
Using the Dynamic Programming Principle (\ref{dpp}) for $W\equiv
V_L[g]$, we can find a sequence $\alpha_n\in {\cal A}$ so that \[
W(\bar x)-W(\phi(t,\bar x,\alpha_n))\le J[g](\bar
x,t,\alpha_n)+\frac{1}{n}.\]   A reapplication of Lemma
\ref{compactness} on the sequential compactness of ${\cal A}^r$,
(\ref{dpp}), and the fact that $W\equiv V_L[\hat g]\equiv
V^r_L[\hat g]\equiv V_L[g] \equiv V^r_L[g]$ on ${\cal D}_o$ then
give $\bar\alpha\in {\cal A}^r$ such that
\begin{eqnarray}
g_o t &\le & \int_0^t \hat g^r(\phi(s,\bar
x,\bar\alpha),\bar\alpha(s))\, {\rm d}s\nonumber\\
&\le& W(\bar x)-W(\phi(t,\bar x,\bar\alpha))\nonumber\\
&\le& \int_0^t  \gamma^r(\phi(s,\bar x,\bar\alpha))\, {\rm
d}s\nonumber\\ &\le & tL_\gamma \max\{||\phi(s,\bar
x,\bar\alpha)-\bar x||: 0\le s\le t\} \; \le\; t\eps,\nonumber
\end{eqnarray}
contradicting the choice of $\eps$.  This shows that $V:=1-e^{-W}$
satisfies the requirement. \er

\br{compare}     Notice that Corollary \ref{prop2} applies to
cases where (\ref{poss}) is not satisfied, and which therefore are
not covered by \cite{CGW01}.   For example, take any bounded ULES
$f$ and
\[g(x,a)\equiv
\twoif{\frac{1}{1+r}\frac{||x||}{r},}{||x||\le r}
{\frac{1}{1+||x||},}{||x||>r}
\] where $r$ is the constant in the ULES definition,
in which case the $g$ quasi-stability of $\R^N$ is easily shown.
In this way, Corollary \ref{prop2} generalizes the uniqueness
results of \cite{CGW01} on (\ref{gz}) by establishing stronger
conclusions that also cover a larger class of cost functions $g$.
Moreover, Remark \ref{rknewly} shows that Corollary \ref{prop2}
gives PDE characterizations for a strictly larger class of
Lyapunov functions than is covered by \cite{CGW01}.\er

\br{twok} Recall that by the ``${\cal KL}$-Lemma'' (cf.
\cite{S98a}), each $\beta_f\in {\cal KL}$ admits
$\alpha_1,\alpha_2\in {\cal K}^\infty$ such that
\begin{equation}
\label{break} \beta_f(r,t)\; \; \le\; \;  \alpha_2(\alpha_1(r)
e^{-t})\; \; \; \; \forall r,t\ge 0
\end{equation}
The results we gave in this note remain true if we replace the
assumption (\ref{rella}) with the condition that
 there be positive constants $\eps$ and
$\delta$ such that
\begin{equation}
\label{liftt} 0\; \; \le g(x,a)\; \; \le\; \;
\delta\alpha^{-1}_2(||x||)\; \; \; \; \; \forall (x,a)\in
B_{\eps}\times A
\end{equation}
and relax the ULES assumption to the requirements that (i) $f$ is
ULAS and (ii) the requirement $(\star\star)$ holds for $\beta_f\in
{\cal KL}$ satisfying (\ref{break})
 (and keep the rest of the assumptions the same).
Condition (\ref{liftt}) is used in \cite{CGW00} to study the Zubov
equation for ULAS dynamics, but under the extra conditions
(\ref{poss}) and $(A_6)$. (In \cite{CGW00}, $f$ is assumed to
satisfy a more general ULAS-like condition, involving a general
 compact
attraction set $D$ instead of the origin. The arguments we gave
above can  be adapted to this slightly more general case.)
 Notice
however that the results we gave on verification functions $w$ in
$\S$\ref{prelim} (for $\ell=-g$) remain  true for general
asymptotically null dynamics, even if (\ref{liftt}) is not assumed
(cf. Remark \ref{con}).\er

\section{Illustrations} \label{illu}

This section illustrates our results by revisiting a simple
example from \cite{CGW01}, where the ULES dynamics is
\begin{equation} \label{lift} f(x,a)=(-x_1+ax^2_1,-x_2+ax^2_2),\;
\; \; \; (x,a)\in \R^2\times [-1,+1]
\end{equation}
(We can take $f\equiv 0$ outside a suitable set containing
$(-1,+1)^2\times A$ to bound the dynamics.)   If we choose
$g(x)=||x||^2$, then the positivity condition (\ref{poss}) on $g$
holds, and \cite{CGW01} shows that
\[
V_L(x)= \twoif {-\ln(1-x_1)-\ln(1-x_2)-x_1-x_2,}{x_1\ge -x_2}
{-\ln(1+x_1)-\ln(1+x_2)+x_1+x_2,}{x_1\le -x_2}
\]
on ${\cal D}={\cal D}_o=(-1,1)^2$.   By \cite{CGW01}, this
function is the unique continuous {\em positive  solution} of
\begin{equation} \label{exeq} x\cdot Dw(x) -
|x^2_1(Dw(x))_1 + x^2_2(Dw(x))_2|-||x||^2\; =\; 0 \end{equation}
on ${\cal D}_o$ that satisfies
\begin{equation}
\label{side} (SC_w),\; \; w \; \text{is\ bounded-from-below,} \;
\text{and}  \; w(x)\to +\infty\; \text{as}\; x\to\partial({\cal
D}_o)\end{equation}
  On the other hand, Theorem
\ref{thm2} proves a stronger result, namely, that $V_L$ is in fact
the unique {\em solution} of (\ref{exeq}) on ${\cal D}_o$ in the
class of functions $w:{\cal D}_o\to\R$ that satisfy (\ref{side}).
In particular, there are no  discontinuous solutions of
(\ref{exeq}) on ${\cal D}_o$ that satisfy (\ref{side}).

Our results also apply to more general cases. For example, change
$g$ to  \begin{equation}\label{gee}\hat g(x,a)\equiv
||x||^2\Psi(x), \; \text{where}\;
\Psi(x)=\sqrt{|x_1-3/4|}+\sqrt{|x_2-3/4|}.\end{equation}
  By Corollary
\ref{prop1}, and the fact that $\Psi(x)\ge
|x_1-\frac{3}{4}|+|x_2-\frac{3}{4}|$ for all $x\in{\cal D}_o$ with
$x_1,x_2\ge 0$, one shows that the Lyapunov function $V_L[\hat g]$
still satisfies (\ref{asy}), and $V_L[\hat g] $ is still a robust
Lyapunov function for $f$ on ${\cal D}_o$ (cf. \cite{BCD97}).
Also, $V_L[\hat g]$ is still a solution of
\begin{equation}
\label{exeqq} x\cdot Dw(x)- |x^2_1(Dw(x))_1 + x^2_2(Dw(x))_2|
-||x||^2\Psi(x) \; =\; 0\end{equation}
 on ${\cal D}_o$, and it is the unique solution of this equation
 on ${\cal D}_o$
 that satisfies
(\ref{side}), by Theorem \ref{thm2}. Since the Dynamic Programming
Principle still holds on $\R^N$,
 $\check V_L[\hat g]$ is the unique bounded
solution of the generalized Zubov equation
\begin{equation}
\label{zex} ||x||^2\Psi(x)\left[w(x)-1\right] + x\cdot Dw(x)-
|x^2_1(Dw(x))_1 + x^2_2(Dw(x))_2|\; =\; 0
\end{equation}
on $\R^N$ that satisfies $(SC_w)$ (cf. \cite{BCD97} and
$\S$\ref{appl3}). Since the Lipschitz and positivity conditions on
$g$ are no longer satisfied, these results do not follow from the
known results. The preceding results remain true if $\Psi$ in
(\ref{gee}) is replaced by
\[\Psi(x)=|x_1-3/4|+|x_2-3/4|.\]
 More
generally, for any  $g:\R^2\times [-1,+1]\to[0,\infty)$ satisfying
the standing hypotheses, we get  a maximal cost type robust
Lyapunov function $V_L$  for (\ref{lift}) and corresponding PDE
characterizations for $V_L$, $\check V_L$, and ${\cal D}_o$.

\section{Conclusion}
\label{concl} This note analyzed a  class of explicit robust
Lyapunov functions of maximal cost type for uniformly locally
asymptotically stable dynamics.  These Lyapunov functions were
shown to be unique viscosity solutions of  first order equations,
subject to appropriate side conditions. The Kru\v{z}kov
transformations of these Lyapunov functions are unique solutions
of the generalized Zubov equation (\ref{gz}) introduced in
\cite{CGW01}. Since we allowed general cost functions $g$, these
uniqueness results do not follow from known results on first-order
viscosity solutions.  The uniqueness characterizations were used
to give sublevel set characterizations for the robust domain of
attraction. As a byproduct, we gave new PDE characterizations for
the variable interest rate infinite horizon minimal cost function
for cases where the Lagrangian could be negative, including the
case where the discount rate is identically zero and the dynamics
is unstable.   This result is of independent interest, because it
allows infinite horizon problems with unbounded cost functions
where cost minimization can take place in one part of the state
space while maximization takes place in the rest of the state
space. One could consider the question of what subset of the set
of all robust Lyapunov functions can be expressed as unique
solutions of Zubov PDE's.  As shown in Remark \ref{rknewly}, our
allowing  general nonnegative $g$ increased the size of this
subset. One could also consider the question of how  the numerical
analysis of the generalized Zubov equation for positive cost
functions $g$ (cf. \cite{CGW00a}) can be extended to the case of
general nonnegative $g$, and in particular, how  allowing
degenerate costs $g$ affects the computation of the sublevel set
${\cal D}_o=\check V^{-1}_L([0,1))$.
  Research on these questions is
ongoing.

\begin{appendix}
\section{Appendix: Quasi-Stability and Existence of
Nearly Optimal Trajectories}


 This appendix proves the following results on quasi-stability
 and nearly-optimal trajectories used in
 $\S\S$\ref{prelim}-\ref{appl}:

\bp{propi}Let $(A_1)$-$(A_4)$ and $(A_6)$ hold and
$\calo\subseteq\R^N$.  Then $\calo$ is relaxed $g$ quasi-stable
iff $\calo$ is $g$ quasi-stable.\ep

\bp{propii}Let $(A_1)$, $(A_2)$,$(A_4)$, and $(A_5)$ hold,  with
$h\equiv 0$ and ${\rm Null}(||f(0,\cdot)||)\ne \emptyset$,
$\cg\subseteq\R^N$ be $f$-invariant, and $w:\cg\to\R$ be a
solution of (\ref{HJBE}) on $\cg\setminus\{0\}$ satisfying
$(SC_w)$.  Let $\bar x\in \cg$, $i\in \N$,  and $\eps>0$ be given,
and define the functions $E_j$ by (\ref{ejay}). Then,
\begin{equation}
\label{zee1} {\cal Z}_i:= \left\{\begin{array}{l}(t,\alpha): \;
0\le t\le 1, \; \alpha\in {\cal A}, \; w_\star(\bar x)\ge \int_0^t
\ell(\phi(s,\bar x,\alpha), \alpha(s))\, {\rm
d}s\\+w_\star(\phi(t,\bar x,\alpha))- E_i(t)\end{array}\right\}
\end{equation}
contains an element of the form $(1, \bar \alpha)$.\ep

We start by proving Proposition \ref{propi}, which  follows from
the following special case of Theorem 1 in \cite{ISW02}:

\bl{lemmai}  Let $A$ be a compact metric space, and let $f$ be
uniformly locally Lipschitz (cf. $(A_2)$).   Let $x\in\R^N$, let
$\alpha\in {\cal A}^r$ be such that $\phi(\cdot,x,\alpha)$ has
domain $[0,\infty)$, and let $r:[0,\infty)\to (0,\infty)$ be
continuous. Then there exist $\beta\in {\cal A}$ and $\eta^o\in
B_{r(0)}(x)$ such that
$||\phi(t,x,\alpha)-\phi(t,\eta^o,\beta)||\le r(t)$ for all
 $t\ge 0$.\el

To prove Proposition \ref{propi}, we adapt the idea from
\cite{ISW02} of putting an approximating trajectory in a tube
around the reference trajectory which has a vanishing radius,
along with an augmentation of the dynamics used in \cite{FR00}.
Extra care is taken  to make sure that not only the reference
trajectory is approximated, but also the {\em integrated cost} of
the trajectory  is approximated.  The details are as follows. Let
$\calo\subseteq\R^N$ be open. If $\calo$ is relaxed $g$
quasi-stable, then it is also $g$ quasi-stable, since as we
remarked above,  ${\cal A}\subseteq {\cal A}^r$. Conversely,
assume $\calo$ is $g$ quasi-stable and $x\in \calo$. Let
$\alpha\in {\cal A}^r$ be such that
$J[g](x,+\infty,\alpha)<+\infty$.  The proposition will follow
once we show that $\phi(t,x,\alpha)\to 0$ as $t\to +\infty$.  We
apply Lemma \ref{lemmai} with the choices
\[
r(t):=\left(\frac{1\wedge \dist(x,\partial\calo)}{2}\right)e^{-t}
\]
and the {\em augmented dynamics}
\begin{equation}
\label{dug1} \twoif{\dot x(t)=f(x(t),a(t))}{}{\dot
y(t)=g(x(t),a(t))}{},\; \; a\in {\cal A}
\end{equation}
Notice that the trajectory for (\ref{dug1}), $a\in {\cal A}^r$,
and any initial position $(\bar x,\bar y)$ has the form
\[
\left(\phi(t,\bar x,a),\, \bar y+\int_0^t g^r(\phi(s,\bar
x,a),a(s))\, {\rm d}s\right)
\]
Using the initial value $(x,0)\in \R^{N+1}$, Lemma \ref{lemmai}
gives an input $\beta\in {\cal A}$ and an initial value
$\eta^o=(\eta^o_L, \eta^o_R)\in\R^{N+1}$ satisfying
\begin{eqnarray}\label{half}
||\phi(t,\eta^o_L,\beta)-\phi(t,x,\alpha)||^2 &+&
\left\vert\eta^o_R+\int_0^tg(\phi(s,\eta^o_L,\beta),\beta(s))\,
{\rm d}s\right.\nonumber\\&-&
\left.\int_0^tg^r(\phi(s,x,\alpha),\alpha(s))\, {\rm d}s\right\vert^2\nonumber\\
&\le& r^2(t)\; \; \forall t\ge 0\end{eqnarray} In particular,
$||\eta^o_L-x||\le \dist(x,\partial\calo)/2$, so  $\eta^o_L\in
\calo$. Since $g$ is nonnegative, (\ref{half}) gives
\[\begin{array}{l}
\int_0^\infty g(\phi(t,\eta^o_L,\beta),\beta(t))\, {\rm d}t \; \;
\le\; \;  |\eta^o_R|+  \int_0^\infty
g^r(\phi(t,x,\alpha),\alpha(t))\, {\rm d}t
 \; <\; \infty\, .\end{array}
\]
Since $\calo$ is $g$ quasi-stable, we get
\[\phi(t,\eta^o_L,\beta)\to 0\; \; \text{as}\; \; t\to +\infty,\]
 so (\ref{half})
gives
\begin{eqnarray}
||\phi(t,x,\alpha)||&\le &
||\phi(t,x,\alpha)-\phi(t,\eta^o_L,\beta)|| +
||\phi(t,\eta^o_L,\beta)|| \nonumber\\ &\le & e^{-t}+
||\phi(t,\eta^o_L,\beta)||\; \to\; 0\; \; \text{as}\; \;
t\to\infty,\nonumber
\end{eqnarray}
as needed.

We turn next to the proof of Proposition \ref{propii}, which is a
generalization  of the proof of the theorem of \cite{MS00}.   We
assume $w$ is continuous and $i=1$, the general case being
similar.  (See Theorem \ref{thm2} for the justification for this
continuity assumption.) Note that ${\cal Z}_1$ is partially
ordered by
\begin{equation} \label{relation} (t_1,\alpha_1)\sim(t_2,
\alpha_2)  \; \; \text{iff} \; \; \left[t_1\le t_2 \; \text{and}\;
\alpha_2\lceil[0,t_1]\equiv \alpha_1\;
\text{a.e.}\right]\end{equation}  Since $\cg$ is $f$-invariant,
the argument from \cite{MS00} shows that every totally ordered
subset of ${\cal Z}_1$ has an upper bound in ${\cal Z}_1$. (If
$\{(t_j,\alpha_j)\}$ is totally ordered in ${\cal Z}_1$, then its
upper bound in ${\cal Z}_1$ is $(\tilde t,\tilde \alpha)$, where
$\tilde t=\sup_j t_j$ and $\tilde\alpha(t):=\alpha_j(t)$ for a.e.
$t\in [0,t_j]$.) It follows from Zorn's Lemma that ${\cal Z}_1$
contains a maximal element $(\bar t, \bar \alpha)$. We will now
show that $\bar t=1$. We can assume that $\phi(\bar t,\bar
x,\bar\alpha)\ne0$ (since $\ell(0,a)\equiv w(0)=0$ and $f(0,\bar
a)= 0$ for some $\bar a\in A$).
  Since $\cg$ is $f$-invariant,
 $\phi(\bar t,\bar x, \bar\alpha)\in \cg$.
 Let $B$ be an open set
containing $\phi(\bar t,\bar x,\bar\alpha)$ whose closure lies in
$\cg\setminus\{0\}$. Suppose that $\bar t<1$, set $q:=\phi(\bar
t,\bar x,\bar\alpha)$, and pick $\delta \in [0, \dist(q,
\partial B)/2)$.  Since $f$ is bounded, it follows that
$T_\delta(q)$, as defined in (\ref{ends}), is positive or
$+\infty$.  By the second part of Lemma \ref{lemma2}, it follows
that there is a $t\in (0,1-\bar t\, )$ and a $\beta\in {\cal A}$
so that
\begin{eqnarray}
\label{newone}
 w(\phi(\bar t,\bar x,\bar\alpha))&\ge& \int_0^t
\ell(\phi(s,\phi(\bar t,\bar x,\bar\alpha), \beta), \beta(s))\,
{\rm d}s\nonumber\\ &+& w(\phi(t, \phi(\bar t,\bar x,\bar\alpha),
\beta))\\\nonumber& - & E_1(\bar t+t)\; +\; E_1(\bar t\, )
\end{eqnarray}
 and so that $\phi(s, \phi(\bar t,\bar x,\bar\alpha),
\beta)\in B$ for all $s\in [0,t]$.   Let $\beta^\sharp$ denote the
concatenation of $\bar \alpha\lceil [0, \bar t\, ]$ followed by
the input $\beta$. If we now combine (\ref{newone}) with the
inequality in (\ref{zee1}) with the choice $\alpha=\bar\alpha$,
then we get
\begin{equation}
\label{thirdnewone}
 w(\bar x)\; \; \ge\; \; \int_0^{\bar
t+t} \ell(\phi(s, \bar x,\beta^\sharp), \beta^\sharp(s))\, {\rm
d}s\; +\; w(\phi(\bar t+t,\bar x, \beta^\sharp))\; -\; E_1(\bar
t+t).
\end{equation}
Since $t$ was chosen so that $\bar t+t<1$, we conclude from
(\ref{thirdnewone}) that \[\left(\bar t+t, \beta^\sharp\right)\in
{\cal Z}_1.\] Since $\beta^\sharp$ is an extension of $\bar
\alpha$, this contradicts the maximality of the pair $(\bar t,
\bar \alpha)$. Therefore, $\bar t=1$.  This proves Proposition
\ref{propii} and completes the proof of Theorem \ref{thm1} for
$h\equiv 0$.  The proof for general $h$ is similar.

\br{det} It should be emphasized that Lemma \ref{lemmai} is {\em
not} an extension of the Filippov-Wa\v{z}ewski Relaxation Theorem,
since the original and approximating trajectories are allowed to
have different starting values. Theorem 1 in \cite{ISW02} is
stated in terms of time-dependent differential inclusions $\dot
x\in F(t,x)$.  To get Lemma \ref{lemmai} from this theorem, first
take $F(t,x)\equiv f(x,A)$, so \[\overline{{\rm co}}
(F(t,x))\equiv f^r(x,A^r).\]  The hypotheses of the theorem are
satisfied since we are assuming that $f$ is locally Lipschitz and
$A$ is bounded. The theorem allows us to approximate any solution
$x$ of $\dot x\in \overline{{\rm co}}(F(t,x))$ by a solution $y$
of $\dot y\in F(t,y)$ in such a way that $||x(t)-y(t)||\le r(t)$
for all $t\ge 0$.  Now take $\alpha\in {\cal A}^r$ and apply the
preceding to the solution $x(t):=\phi(t,x,\alpha)$ of $\dot x\in
\overline{{\rm co}}(F(t,x))$. Since the sets $A(t):=\{a\in A:
f(y(t),a)=\dot y(t)\}$ are measurable, Filippov's Lemma  gives an
input $\beta\in {\cal A}$ such that
$y(\cdot)=\phi(\cdot,\eta^o,\beta)$, where $\eta_o\in
B_{r(0)}(x)$, which gives Lemma \ref{lemmai}.\er

\section{Appendix: Nonglobal Solutions of Zubov Equation}
In $\S$\ref{appl}, we gave global uniqueness characterizations for
solutions of the PDE (\ref{gz}) on $\R^N$.  As pointed out in
\cite{CGW01}, it can be inconvenient from a practical point of
view to verify that a function is a PDE solution on all of $\R^N$.
This motivates the problem of finding {\em non}global uniqueness
characterizations for solutions of (\ref{gz}), on general open
sets $\calo\subseteq\R^N$. For simplicity, we will assume ${\cal
D}_o\subseteq \calo$.  The following nonglobal uniqueness
characterization extends the results of \cite{CGW01}, $\S$3, to
general continuous $g$. The proof is a localization of the
argument  used to prove  Theorem \ref{thm1}. Recall that if
$F:\R^N\times \R\times\R^N\to\R$ is continuous, then a locally
bounded function $w:\calo\to\R$ is called a {\bf (viscosity)
supsersolution} (resp., {\bf subsolution}) of
\[F(x,w(x),Dw(x))=0\] on $\calo$ provided that $(C_1)$ (resp.,
$(C_2)$) of Definition \ref{definitionvis} is satisfied. We will
use the fact (cf. \cite{BCD97}) that if  $(A_1)$-$(A_5)$ hold,
then any subsolution (resp., supersolution) of (\ref{HJBE}) on a
bounded open set $B$ satisfies conclusion (a) (resp., (b)) of
Lemma \ref{lemma2}.

\bp{lastone} Let $(A_1)$-$(A_4)$ and $(A_6)$-$(A_8)$ hold with $f$
ULES.  Assume that $\calo\subseteq\R^N$ is open, ${\cal
D}_o\subseteq\calo$, and $w:{\rm cl}(\calo)\to\R$ is bounded. \bi
\item[(i)]
If $w$ is an upper semicontinous function which is a  subsolution
of (\ref{gz}) on $\calo$ satisfying $w(0)\le 0$ and $w\equiv 1$ on
$\partial (\calo)$, then $w\le \check V_L[g]$ on $\calo$.
\item[(ii)]
If $w$ is a lower semicontinous function which is a supersolution
of (\ref{gz}) on $\calo$ satisfying $w(0)\ge 0$ and $w\ge 1$ on
$\partial (\calo)$, then $w\ge \check V_L[g]$ on $\calo$. \ei\ep

\bpr (i)  By the proof of Theorem \ref{thm1} (with $\ell\equiv -g$
and $h\equiv g$), $w\le \check V_L[g]$ on ${\cal D}_o$. Let $\bar
x\in \calo\setminus {\cal D}_o$.  It remains to show that \[w(\bar
x)\le \cvlg(\bar x).\]  Assume the contrary and pick $\eps>0$ such
that
\begin{equation}\label{contra}
w(\bar x)\ge \cvlg(\bar x)+\eps
\end{equation}
Define the functions $E_i$ by (\ref{ejay}) for $i\in \N$.
  Set $\bar\phi(0)=\bar x$, and
\[
J_w(x,t,\alpha) =G(x,t,\alpha)-1-G(x,t,\alpha) w(\phi(t,x,\alpha))
\]
wherever the RHS is defined, where $G$ is as defined by
(\ref{DP}). We will inductively define the sets
\[
{\cal Z}_{\calo, i+1}:= \left\{\begin{array}{l} (t,\alpha)\in
[0,1]\times {\cal A}:\; \; \phi(s,\bar\phi(i),\alpha)\in \calo \;
\, \forall s\in [0,t],\\-w(x)\; \ge\;  J_w(\bar \phi(i),t,\alpha)
-\frac{1}{2} E_{i+1}(t)\end{array}\right\}
\]
for all $i\in \N_o$.  The set ${\cal Z}_{\calo,1}$ is partially
ordered by $\sim$ as defined in (\ref{relation}).  Let \[\hat
T:=\{(t_j,\alpha_j)\}\] be a totally ordered subset of ${\cal
Z}_{\calo,1}$.  If $\phi(t_j,\bar x,\alpha_j)$ converges to  point
in $\partial(\calo)$, then the semicontinuity of $w$ implies that
$-\cvlg(\bar x)-\eps$ majorizes
\begin{eqnarray}
\label{argu} -w(\bar x) &\ge & G(\bar x,t_j,\alpha_j)-1-G(\bar
x,t_j,\alpha_j)\, w(\phi(t_j,\bar
x,\alpha_j))-\frac{\eps}{2}\nonumber\\
&\ge &G(\bar x,t_j,\alpha_j)-1-G(\bar
x,t_j,\alpha_j)\left(1+\frac{\eps}{4}\right)
-\frac{\eps}{2}\nonumber\\
&\ge &-1-\frac{3\eps}{4}
\end{eqnarray}
for large $j$, so $\cvlg(\bar x)\le 1-\eps/4$. Therefore, $\bar
x\in {\rm dom}(V_L[g])={\cal D}_o$.  This contradicts the choice
of $\bar x$.  By the proof of Proposition \ref{propii}, it follows
that $\hat T$ has an upper bound in ${\cal Z}_{\calo,1}$ (cf.
\cite{MS00}). By Zorn's Lemma, ${\cal Z}_{\calo,1}$ contains a
maximal element $(\bar t,\bar\alpha)$, and the proof of Theorem
\ref{thm1} (with $\ell=-g$ and $g\equiv h$, applied to the
supersolution $-w$ of the Hamilton-Jacobi equation (\ref{HJBE}))
shows that $\bar t=1$. Now set
\[\bar\phi(1)=\phi(1,\bar x,\bar\alpha).\] The preceding argument
gives a maximal element $(1,\bar\alpha_2)\in {\cal Z}_{\calo,2}$.
(To show that each totally ordered set of $(s_j,\alpha_{j,2})$'s
in ${\cal Z}_{\calo,2}$ has an upper bound in ${\cal
Z}_{\calo,2}$, notice that if
$\phi(s_j,\bar\phi(1),\alpha_{j,2})\to p\in
\partial({\calo})$ for some $p\in\R^N$,
then (\ref{argu}) with $t_j$ and $\alpha_j$ replaced by $1+s_j$
and the concatenation of $\bar\alpha$ followed by $\alpha_{j,2}$,
respectively, gives the same contradiction as before.)  Set
\[\bar\phi(2)=\phi(1,\bar\phi(1),\bar\alpha_2).\] This procedure is
iterated exactly as in the proof of Theorem \ref{thm1} and gives
an input $\tilde\alpha\in {\cal A}$ such that
\[w(\bar x)\; \; \le\; \;
1-G(\bar x,M,\tilde\alpha)+G(\bar x,M,\tilde\alpha)\,
w(\phi(M,\bar x,\tilde\alpha))+\frac{\eps}{2}\; \; \forall
M\in\N\] and such that $\phi(t,\bar x,\hat\alpha)\in \calo$ for
all $t\ge 0$.  A reapplication of (\ref{lom}) gives \[w(\bar x)\;
\le\; 1-G(\bar x,+\infty,\tilde \alpha)+\frac{\eps}{2}\; \le\;
 \cvlg(\bar x)+\frac{\eps}{2},\] which again contradicts the
choice of $\eps>0$ in (\ref{contra}).  Therefore, \[w(\bar x)\le
\cvlg(\bar x),\] which proves  (i).

\noindent (ii) Let $\bar x\in \calo$.  For any $\alpha\in {\cal
A}$, define the exit times
\[
\tilde t(x,\alpha):=\inf\{t\ge 0: \phi(t,x,\alpha)\not\in
\calo\}\in [0,+\infty]
\]
A repeated application of the first part of Lemma \ref{lemma2}
with $\ell=-g$ (cf. \cite{M00}) gives
\begin{equation}
\label{tiltil}
 w(\bar x)\; \; \ge\; \;
1-G(\bar x, t,\alpha)+G(\bar x,t,\alpha)\, w(\phi(t,\bar
x,\alpha))
 \; \; \forall \alpha\in {\cal
 A}
\end{equation}
for all finite $t\in (0,\tilde t(\bar x,\alpha)]$.
 If $\tilde t(\bar x,\alpha)$ is finite for some
$\alpha\in {\cal A}$, then (\ref{tiltil}) with the choice
$t=\tilde t(\bar x,\alpha)$ gives \[w(\bar x)\ge 1\ge \cvlg(\bar
x).\] Otherwise, (\ref{lom}) and a passage to the liminf as $t\to
+\infty$ in (\ref{tiltil}) for fixed $\alpha$ gives \[w(\bar x)\ge
1-G(\bar x,+\infty,\alpha)\]  for all $\alpha\in {\cal A}$, so
$w(\bar x)\ge \cvlg(\bar x)$, as needed.
 \epr
\br{lastlast} If we put (i)-(ii) of the previous proposition
together, then we get nonglobal PDE characterizations for $\cvlg$
which extend the results of \cite{CGW01}.\er
\end{appendix}

\noindent{\bf Acknowledgements.} Part of this work was carried out
in 1999, while the author was a Predoctoral Fellow at Rutgers
University in New Brunswick, NJ. The author thanks Eduardo D.
Sontag and H\'ector J. Sussmann for helpful discussions during
this period.  The author was supported in part by USAF Grant
F49620-98-1-0242 and DIMACS Grant NSF CCR91-19999.  The author
thanks the three anonymous referees for their careful reading of
this manuscript and for their extensive comments.

\end{document}